\documentclass[twoside,11pt]{article}

\usepackage[preprint]{jmlr2e}

\usepackage[utf8]{inputenc}

\usepackage{amsmath,amssymb,amsfonts,mathtools}
\usepackage{subcaption}
\usepackage{algorithm}
\usepackage{algpseudocode}
\usepackage{bm}
\usepackage{makecell}
\usepackage{booktabs}
\usepackage{scalerel}
\usepackage{dsfont}
\usepackage{multicol}
\usepackage{multirow}
\usepackage[table]{xcolor}
\usepackage{enumitem}
\usepackage{pgfplots}
\usepackage{float}
\newcommand{\E}{\mathbb{E}}
\newcommand{\Var}{\mathrm{Var}}
\newcommand{\Prob}{\mathbb{P}}

\ShortHeadings{Nonparametric Estimation of Joint Entropy}{Ho, Park and Oh}
\firstpageno{1}

\title{Nonparametric Estimation of Joint Entropy\\
via Partitioned Sample-Spacing}

\author{%
  \name Jungwoo Ho \email hoo0321@yonsei.ac.kr \\
  \addr Department of Statistics and Data Science, Yonsei University
  \AND
  \name Sangun Park \\
  \addr Department of Statistics and Data Science, Yonsei University
  \AND
  \name Soyeong Oh \\
  \addr Department of Statistics and Probability, Michigan State University
}

\begin{document}

\maketitle

\begin{abstract}
We propose a nonparametric estimator of multivariate joint entropy based on partitioned sample spacing (PSS). The method extends univariate spacing ideas to $\mathbb{R}^d$ by partitioning into localized cells and aggregating within-cell statistics, with strong consistency guarantees under mild conditions. In benchmarks across diverse distributions, PSS consistently outperforms $k$-nearest neighbor estimators and achieves accuracy competitive with recent normalizing flow–based methods, while requiring no training or auxiliary density modeling. The estimator scales favorably in moderately high dimensions ($d=10$–$40$) and shows particular robustness to correlated or skewed distributions. These properties position PSS as a practical and reliable alternative to both
$k$NN and NF-based entropy estimators, with broad utility in
information-theoretic machine learning tasks such as total-correlation
estimation, representation learning, and feature selection.

\end{abstract}


% ----------------- BODY -----------------

\section{Introduction}

Entropy and mutual information are fundamental quantities in data analysis and machine learning, as they quantify uncertainty and dependence between variables \citep{Shannon1948}. Unlike correlation, which captures only linear relationships, they can detect general nonlinear dependencies and remain invariant under monotone transformations. These properties have made entropy and mutual information central tools for feature selection, representation learning, independent component analysis and causal discovery.  

Estimating entropy and mutual information from samples is challenging. Among nonparametric approaches, the $k$-nearest neighbor (kNN) estimator is widely used \citep{kraskov2004}. It requires no distributional assumptions and adapts naturally to new data, but its accuracy degrades in multivariate settings, particularly under strong correlation or high dimension. More broadly, there is a long line of work on nonparametric estimation of divergences and conditional information using $k$NN statistics \citep{poczos}, which provides consistency guarantees but suffers from the same scalability issues. Recent refinements such as bias-corrected weighting \citep{BerrettSamworthYuan2019} improve asymptotic properties, yet all kNN variants rely on costly neighbor searches and remain sensitive to dependence.  

Another line of work leverages copula decompositions. Ariel and Louzoun~\citep{copula} proposed the Copula Decomposition Entropy Estimator (CADEE), which separates marginal and copula contributions using Sklar’s theorem. By recursively estimating copula entropy on the compact $[0,1]^d$ domain, CADEE achieves strong robustness and scales better than kNN, particularly when supports are mixed or unbounded.  

In parallel, trainable normalizing flow methods have emerged. \citet{AoLi2022} introduced an estimator that learns a bijection to a tractable base distribution, while \citet{ButakovTolmachev2024} proposed MIENF, which fits a pair of flows with finite-sample guarantees. \citet{saad} developed inference-based entropy estimators (EEVI) that compute bounds within probabilistic generative models. 
While such estimators can reduce bias, they require significant computational resources for training or access to tractable model densities, and fail when the model family is misspecified.  

In this work, we revisit Vasicek’s classical sample-spacing idea \citep{Vasicek}\citep{Dudewicz1986161} and extend it to the multivariate setting by applying it locally within partitions of the sample space. The resulting estimator is fully nonparametric, avoids both neighbor search and model training, and remains robust under heavy tails and nonlinear dependence. Under mild assumptions—boundedness and continuity of the density within each partition—we prove strong consistency, with an additional moment condition for entropy. Empirically, our method consistently outperforms kNN in mean-squared error, improves upon CADEE in correlated and moderately high-dimensional settings, and achieves accuracy competitive with normalizing-flow approaches while being substantially more efficient. These results position partitioned sample spacings as a practical and
data-efficient tool for modern information-theoretic learning tasks.  

\section{Preliminary}

\subsection{Nearest Neighbor--Based Estimators}

One widely used class of nonparametric entropy estimators relies on $k$-nearest neighbor (kNN) distances.  
Given $n$ i.i.d.\ samples $X_1,\dots,X_n \in \mathbb{R}^d$, let $R_{i,k}$ denote the distance from $X_i$ to its $k$th nearest neighbor.  
These distances provide local volume estimates, which can be inverted to obtain density estimates and hence entropy.

\paragraph{Kozachenko--Leonenko (KL) estimator.}
The classical estimator of \citet{Kozachenko} assumes local uniformity within a ball of radius $R_{i,k}$ around $X_i$.  
With $V_d = \pi^{d/2}/\Gamma(1+d/2)$ the volume of the $d$-dimensional unit ball, the KL estimator is
\[
\hat H_{\mathrm{KL}}
= -\psi(k) + \psi(n) + \log V_d
+ \frac{d}{n}\sum_{i=1}^n \log R_{i,k},
\]
where $\psi$ is the digamma function.  
Despite its simplicity, the KL estimator can suffer from large bias in high dimensions \citep{gao2015efficientestimationmutualinformation}.

\paragraph{Kraskov--St\"ogbauer--Grassberger (KSG) estimator.}
\citet{kraskov2004} proposed a modification using hyper-rectangles instead of balls.  
For each $X_i$, let $R_{i,j}$ be the distance to its $k$th nearest neighbor along coordinate $j$.  
Defining a rectangle with side lengths $\{R_{i,j}\}$, the KSG entropy estimate becomes
\[
\hat H_{\mathrm{KSG}}
= -\psi(k) + \psi(n) + \frac{d-1}{k}
+ \frac{1}{n}\sum_{i=1}^n \sum_{j=1}^d \log R_{i,j}.
\]
Originally introduced for mutual information, KSG often improves stability compared to KL, especially under strong dependence.

\paragraph{Truncated estimators and Normalizing Flows.}
Boundary bias can be severe when the support is compact.  
To address this, \citet{AoLi2022} introduced truncated versions (tKL, tKSG), where neighbor cells are clipped to the support $[0,1]^d$.  
They proved that tKL is unbiased for uniform distributions.  
Building on this property, they proposed a normalizing flow transformation that maps samples approximately to uniform, and then applies the truncated estimator with a Jacobian correction.  
This yields the so-called \emph{uniformizing mapping} (UM) estimators (UM-tKL, UM-tKSG), which achieve substantially lower bias in high dimensions.  
Normalizing flows themselves are neural network models of invertible transformations \citep{RezendeMohamed2015, Papamakarios2021}.  
They require GPU-based training with many gradient updates, making them far more computationally demanding than nonparametric methods.  
Thus NF-based entropy estimation can reach high accuracy but at the cost of significant training overhead.

\subsection{Univariate Case}
Let $x_{(1)}<x_{(2)}<\cdots<x_{(n)}$ be order statistics of an i.i.d. sample  $x_1,x_2,\dots,x_n$ drawn from a continuous random variable $X$ with probability density function $f(x)$. The difference $x_{(i+m)}-x_{(i)}$ is called m-spacing \citep{Learned-Miller20041271}. Using this idea, \citet{park_and_park} proposed the following density estimator based on sample spacings:
\begin{equation*}
\hat{f}_n(x) = 
\begin{cases}
\displaystyle\frac{2m}{n\big(x_{(i+m)} - x_{(i-m)}\big)} 
  & \text{for } \xi_i < x \le \xi_{i+1}\\
0 & \text{otherwise}
\end{cases}
\end{equation*}
 where $\xi_{i}=\frac{x_{(i-m)}+\cdots+x_{(i+m-1)}}{2m}$, and $x_{(i)}=x_{(n)}$ for $i>n$. Note that for $i<1$, $x_{(i)}=x_{(1)}$ and for $i>n$, $x_{(i)}=x_{(n)}$. Based on this m-spacing method differential entropy is derived as (Vasicek, 1976)
\begin{equation*}
H_v(n, m) = \frac{1}{n} \sum_{i=1}^n \log\left( \frac{n}{2m} (x_{(i+m)} - x_{(i-m)}) \right) 
\end{equation*}
 This entropy estimator converges to $H(f)$ as $n,m \to \infty$ and $\frac{m}{n} \to 0$ \citep{Dudewicz1986161}.

While \citet{Vasicek} proved the consistency of the entropy estimator $H_v(n, m)$, he did not establish the consistency of the corresponding density estimator. We fill this gap by establishing the almost sure consistency of $\hat{f_n}(x)$, under mild assumptions on the spacing parameter $m$.
The key step is to prove that the averaged grid point $\xi_i$ converges to the order statistic $X_{(i)}$, ensuring that the estimator localizes correctly around the sample point. 

\begin{lemma}[Consistency of $\xi_{i}$]
     If $m \rightarrow \infty$ and $m/n \rightarrow 0$ as $n \rightarrow \infty$, then $\xi_{i}-X_{(i)}$ converges to 0 almost surely.
\end{lemma}

See Appendix A for the proof. The choice of the spacing parameter $m$ is critical for balancing bias and variance. To ensure consistency, $m$ must grow with the sample size $n$ \citep{Bert-vanEs}. For our theoretical results and empirical studies, we adopt the widely used rate $m=\lfloor\sqrt{n}+1/2\rfloor$ \citep{Crzc}, which is known to be effective for smooth distributions.

\begin{theorem}[Consistency of $\hat{f}_n(x)$]
If the spacing parameter is set as $m=\lfloor\sqrt{n}+1/2\rfloor$, then the sample-spacing density estimator $\hat{f}_n(x)$ converges almost surely to $f(x)$:
\[
\hat{f}_n(x) \xrightarrow{a.s.} f(x),
\]
where $\xrightarrow{a.s.}$ denotes almost sure convergence.
\end{theorem}

See Appendix B for the proof. This result formally establishes the almost sure consistency of the univariate sample-spacing density estimator, providing a theoretical foundation for its extension to the multivariate and partitioned settings developed in the next sections.

\section{Proposed Extension}

\subsection{Bivariate Case}

Let $(X_1, Y_1), \dots, (X_n, Y_n)$ be i.i.d.\ samples from a continuous joint density $f(x,y)$ defined on $\mathbb{R}^2$. The joint (differential) entropy is
\begin{equation*}
    H(X,Y) = - \iint_{\mathbb{R}^2} f(x,y) \log f(x,y) \, dx \, dy ,
\end{equation*}
which measures the overall uncertainty in the bivariate distribution.

Vasicek’s sample-spacing method was originally proposed for univariate entropy estimation. However, the same idea can be extended to higher dimensions by constructing a piecewise-constant joint density estimator on a grid formed by the order statistics of each variable. This generalization retains the nonparametric nature of the method while enabling it to capture relationships between multiple variables.

Let $x_{(1)} \le \dots \le x_{(n)}$ and $y_{(1)} \le \dots \le y_{(n)}$ be the order statistics of $\{X_i\}$ and $\{Y_i\}$, respectively. Define the averaged grid points $\xi_i = \frac{x_{(i-m_x)} + \cdots + x_{(i+m_x-1)}}{2m_x}$, $\eta_j = \frac{y_{(j-m_y)} + \cdots + y_{(j+m_y-1)}}{2m_y}$, with $x_{(i)} = x_{(1)}$ for $i<1$, $x_{(i)} = x_{(n)}$ for $i>n$, and analogously for $y_{(j)}$.

If $X$ and $Y$ are independent, the joint density factorizes as $f(x,y) = f_X(x) f_Y(y)$. In this case, the joint density estimator can be constructed as the product of two univariate spacing-based estimators $\hat{f}_{n}^{\mathrm{ind}}(x,y)$:
\[
\hat f_{n}^{\mathrm{ind}}(x,y) =
\frac{4 m_x m_y}{n^2 \,(x_{(i+m_x)} - x_{(i-m_x)}) \, (y_{(j+m_y)} - y_{(j-m_y)})}
\]
for $\xi_i < x \le \xi_{i+1}, \ \eta_j < y \le \eta_{j+1}$, $i,j = 1, \dots, n$.  
The corresponding joint entropy estimator is then
\begin{equation*}
    H_{v}^{\mathrm{ind}}(n, m_x, m_y) =
    -\frac{1}{n^2} \sum_{i=1}^n \sum_{j=1}^n
    \log \hat{f}_{n}^{\mathrm{ind}}(x,y).
\end{equation*}

A direct extension of the univariate sample-spacing method to the bivariate case constructs a grid from marginal order statistics. 
This is simple but flawed: marginal grids ignore joint geometry, so dense regions may be covered by large low-density cells and vice versa. 
For correlated data, the estimator underestimates along dense diagonals and overestimates in sparse regions, because marginal orderings distort the grid. 

\begin{figure}[h]
    \centering
    \subfloat[\centering Direct Method ($\ell=1$)]{{\includegraphics[width=6cm]{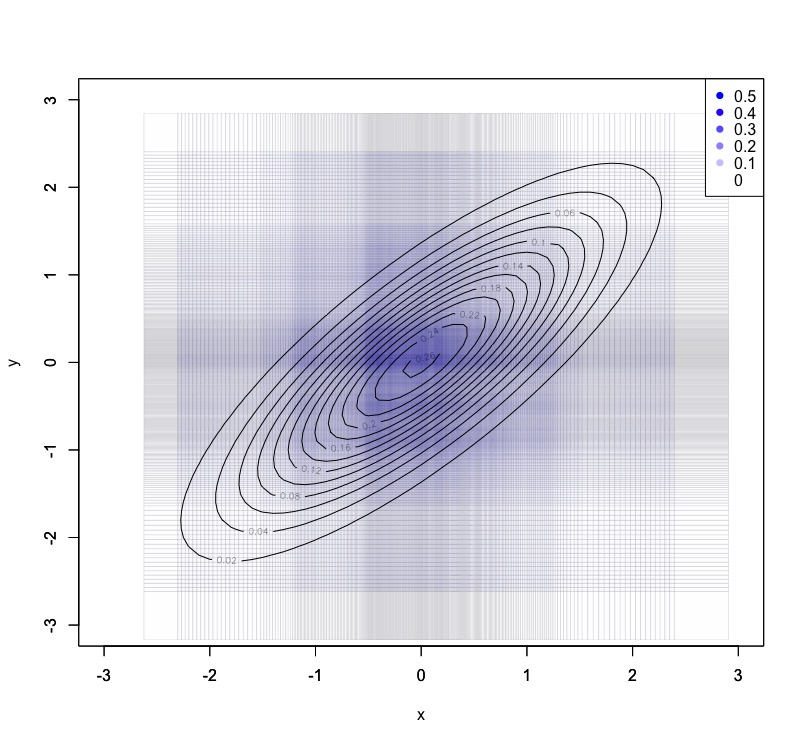}}}%
    \quad \quad
    \subfloat[\centering Partitioned Method ($\ell=3$)]{{\includegraphics[width=6cm]{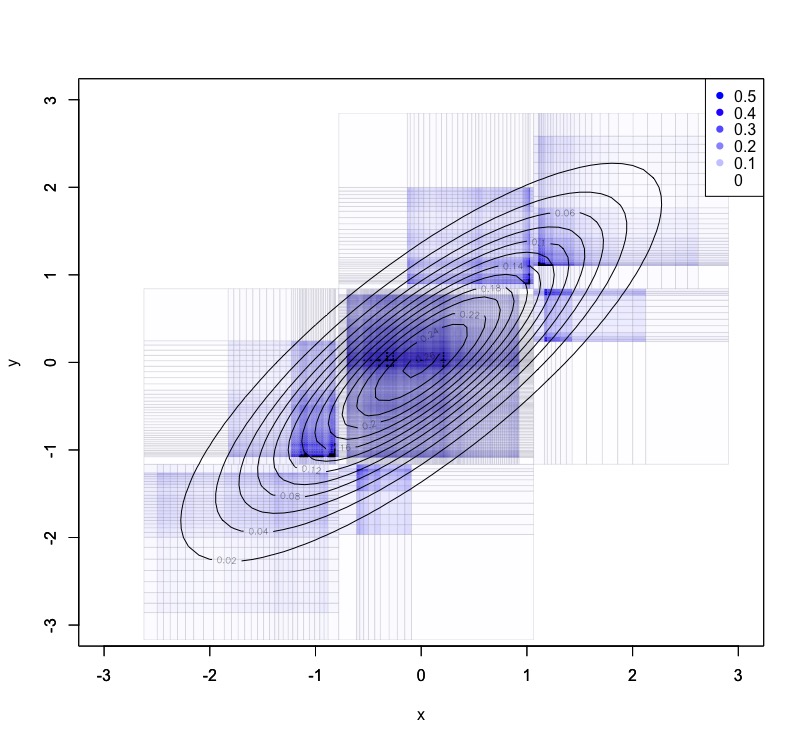} }}%
    \caption{Comparison of density estimation for a bivariate normal with $\rho=0.8$. The direct method (left) fails to capture the correlation, while the partitioned method (right) successfully reflects the joint structure.}
    \label{fig:l1_vs_l3}
\end{figure}

Figure~\ref{fig:l1_vs_l3}(a) shows this failure for a bivariate normal with $\rho=0.8$. 
The direct method does not reflect the linear dependence and even assigns spurious high density to off-diagonal regions. 
To address this, we partition each axis into $\ell$ equal-width segments and apply the sample-spacing estimator locally. 
This prevents mis-assignment across quadrants and better captures the correlation, as seen in Figure~\ref{fig:l1_vs_l3}(b) with $\ell=3$. 

Unlike adaptive partitioning approaches based on independence tests \citep{Darbellay19991315}, or sample-spacing methods assuming radial symmetry \citep{Lee20102208}, 
our method uses fixed, assumption-free partitions, making it both simple and broadly applicable.

\subsection{The Partitioned Density Estimator}

Let $P_k$ be the $k$-th partition, $m_{k}$ be the spacing parameter in $P_k$, $n_{k}$ be the number of samples in that partition, and $(x_{(1)}^k,\dots,x_{(n_k)}^k)$, $(y_{(1)}^k,\dots,y_{(n_k)}^k)$ be the corresponding order statistics for $k=1,\cdots,\ell^2$. To ensure the estimator integrates to 1, we weigh the local density estimate in the $k$-th partition by the sample proportion $\frac{n_k}{n}$. The partitioned sample-spacing density estimator $\hat{f}_{n,\ell}(x,y)$ is:
\begin{equation*}
    \hat f_{n,\ell}(x,y) =
    \frac{n_{k}}{n}\frac{4m_{k}^2}{n_{k}^2(x_{(i+m_{k})}^k-x_{(i-m_{k})}^k)(y_{(j+m_{k})}^k-y_{(j-m_{k})}^k)} 
\end{equation*}
where $(x,y)\in (\xi_{i}^k, \xi_{i+1}^k]\times(\eta_{j}^k,\eta_{j+1}^k]$. The grid points are defined as $\xi_i^k=\frac{x_{(i-m_k)}^k+\dots+x_{(i+m_k-1)}^k}{2m_{k}}$ and $\eta_{j}^k=\frac{y_{(j-m_k)}^k+\dots+y_{(j+m_k-1)}^k}{2m_k}$, with boundary points $\xi_0^k=x_{(1)}^k$ and $\xi_{n_k+1}^k=x_{(n_k)}^k$. Analogously, $\eta_0^k=y_{(1)}^k$ and $\eta_{n_k+1}^k=y_{(n_k)}^k$. The partitioned density estimator is defined locally within each partition.  

To make the estimator more adaptive and reduce estimation bias, the sub-grid is constructed on the empirical support of the data within each partition, i.e.,  on the partition $[\xi_{0}^k,\xi_{n_k+1}^k]\times [\eta_{0}^k,\eta_{n_k+1}^k]$. While this means the density is implicitly zero in the regions between the data's range and the partition boundaries for finite samples, this choice is empirically superior, yielding a significantly lower MSE. As a result of this data-driven construction, the estimator is not guaranteed to integrate to exactly 1 for finite samples; however, it asymptotically integrates to 1, as the empty regions between the data's range and the partition boundaries vanish as $n_k\rightarrow\infty$. 

\begin{proposition}
The partitioned sample-spacing estimator $\hat{f}_{n,\ell}(x,y)$ is asymptotically a valid probability density function, satisfying
\[
\lim_{n \to \infty} \int_{\mathbb{R}^2} \hat{f}_{n,\ell}(x,y)\,dx\,dy = 1.
\]
\end{proposition}

See Appendix C for the proof. This estimator is a valid probability density function and is strongly consistent. We formalize this in the following lemma and theorem. The proof strategy is to first establish, in Lemma~\ref{lem:spacing_consistency}, the almost sure convergence of the univariate $m$-spacing estimators for the conditional distributions within each partition. This result serves as the core ingredient for Theorem~\ref{thm:bivar-consistency}, which combines these components to prove the consistency of the overall joint density estimator.

The assumptions in the lemma are critical for this two-stage argument. The conditions on $\ell(n)$ and $n_k$ ensure that as the total sample size $n$ grows, the partitions shrink appropriately while still accumulating enough data for asymptotic analysis. Similarly, the condition on the spacing parameter $m_k$ is crucial. As established in the preceding section, we adopt the specific rate $m_k = \lfloor\sqrt{n_k} + 1/2\rfloor$, which is required to prove convergence. Hence, we assume the conditions on $\ell(n)$, $n_k$, and the specific rate for $m_k$ hold for all subsequent asymptotic results.

\begin{lemma}\label{lem:spacing_consistency}
Let the event
$C_k \;=\;\bigl\{(X,Y)\in P_k\bigr\}$ and $n_k :=\sum_{v=1}^n \mathbf1\{(X_v,Y_v)\in C_k\}.$ Assume that as $n\to\infty$, $\ell(n)\to\infty,\ \ell(n)^2 = o(n),\
n_k\rightarrow \infty.$ Further, assume the spacing parameter $m_k = \lfloor\sqrt{n_k} + 1/2\rfloor$. For $x\in(\xi_j^k,\xi_{j+1}^k]
,\
y\in(\eta_i^k,\eta_{i+1}^k]$
define the univariate $m$‐spacing estimators
\begin{align*}
\widehat f_{X,k}(x)
&=\frac{2\,m_k}{\,n_k\bigl(x_{k,(i+m_k)} - x_{k,(i-m_k)}\bigr)},\\
\widehat f_{Y,k}(y)
&=\frac{2\,m_k}{\,n_k\bigl(y_{k,(j+m_k)} - y_{k,(j-m_k)}\bigr)}.
\end{align*}
Then for any fixed $(x,y)\in(\eta_i^k,\eta_{i+1}^k]\times(\xi_j^k,\xi_{j+1}^k]$,
\[
\widehat f_{X,k}(x)\;\xrightarrow{a.s.}\;f_{X\mid C_k}(x),
\quad
\widehat f_{Y,k}(y)\;\xrightarrow{a.s.}\;f_{Y\mid C_k}(y).
\]
\end{lemma}
See Appendix D for the proof.
\begin{theorem}[Consistency of $\hat{f}_{n,\ell}(x,y)$]\label{thm:bivar-consistency}
The partitioned bivariate density estimator 
\[
 \hat{f}_{n,\ell}(x,y) = \frac{n_k}{n} \cdot \widehat f_{X,k}(x) \cdot \widehat f_{Y,k}(y),
\]
converges almost surely to $f(x,y)$.
\end{theorem}

See Appendix E for the proof. The almost sure convergence in Theorem~\ref{thm:bivar-consistency} provides the foundation for $L^1$ convergence via Scheffé's theorem \citep{Williams}: if a sequence of nonnegative integrable functions
$g_n$ satisfies (i) $g_n\to g$ a.e. and (ii) $\int g_n \to \int g$, then
$\int |g_n-g|\to 0$. In our setting, Theorem~\ref{thm:bivar-consistency} yields (i) for $\hat f_{n,\ell}$, and
Proposition~1 establishes (ii) with $\int f=1$.

\begin{corollary}[{$L^1$ convergence of $\hat f_{n,\ell}$}]
By Scheffé's theorem, the $L^1$ distance between the estimator and the true density converges to zero:
\[
  \int_{\mathbb{R}^2}\big|\hat f_{n,\ell}(x,y)-f(x,y)\big|\,dx\,dy
  \;\longrightarrow\;0 .
\]
\end{corollary}

\subsection{The Joint Entropy Estimator}

Given the consistent density estimator $\hat{f}_{n,\ell}(x,y)$, a natural and powerful way to estimate the joint entropy is via a ``plug-in'' approach. We replace the true density $f$ with our estimator $\hat{f}_{n,\ell}$ in the definition of entropy and approximate the expectation using the empirical distribution of the $n$ observed samples. This leads to the final entropy estimator:
\begin{equation*}
    \hat{H}_{n,\ell}(X,Y) = -\frac{1}{n}\sum_{v=1}^{n} \log \hat{f}_{n,\ell}(X_v,Y_v).
\end{equation*}
Applying the partition parameter $\ell$ offers several advantages. Since this estimator is computed locally, it adapts well to highly correlated data and effectively captures distributional details such as skewness and kurtosis.

Additionally, the partitioned model significantly improves computational efficiency by reducing the number of effective grid cells. For example, with $n$ samples, the global (unpartitioned) method implicitly forms a grid of size $n^2$, whereas the partitioned model constructs independent grids summing to $\sum_{k=1}^\ell n_k^2$. Since $\sum n_k^2 \ll (\sum n_k)^2$, the computational burden decreases as $\ell$ increases. However, the choice of $\ell$ must also balance the trade-off involving variance, correlation, and sample size.

Notably, the optimal $\ell$ is expected to be invariant to the scale of the data. Because the length of each partition is set as proportional to the range (e.g., $\frac{\max(x) - \min(x)}{\ell}$), the number of samples $n_k$ falling into each partition remains constant under affine transformations. We formally demonstrate this property below.

\begin{proposition}
    The partition parameter $\ell$ is location invariant and scale equivariant.
\end{proposition}

See Appendix F for the proof. Building upon the consistency of the density estimator, we now both establish the almost sure and the $L^1$ convergence of the plug-in entropy estimator. The crucial step is to prove the uniform integrability of the log-density estimator sequence. This requires the globally adopted rate for the spacing parameter, $m_k = \lfloor\sqrt{n_k} + 1/2\rfloor$, as well as a mild integrability condition on the true density $f$ itself. The following theorem formalizes this main result.

\begin{theorem}
Assume that the true density $f(x,y)$ satisfies $\log f \in L^{1+\delta}$ for some $\delta \in (0,1)$. Let the spacing parameter be $m_k=\lfloor\sqrt{n_k}+1/2\rfloor$ with $n_k> 1$, and let the number of partitions satisfy $\ell(n)^2 = o(n)$. Then, as $n \to \infty$, the entropy estimator $\hat{H}_{n,\ell}(X,Y)$ converges to the true entropy $H(f) := \E[-\log f(X,Y)]$ both almost surely and in $L^1$:
\[
\hat{H}_{n,\ell}(X,Y) \xrightarrow{a.s.} H(f) \quad \text{and} \quad\hat{H}_{n,\ell}(X,Y) \xrightarrow{L^1} H(f).
\]
\end{theorem}
See Appendix G for the proof.
\section{Multivariate Extension}
\subsection{Multivariate Density Estimator}

The partitioned sample-spacing method, initially developed for the bivariate case, can be extended to multivariate scenarios. Let \(\mathbf{X} = (X_1, X_2, \dots, X_d)\) represent a \(d\)-dimensional random vector, with \(n\) samples \(\{(x_{v1}, x_{v2}, \dots, x_{vd})\}_{v=1}^n\). The partition parameter \(\ell\) divides each dimension into \(\ell\) segments, resulting in \(\ell^d\) hyper-rectangular partitions, where each partition has a size proportional to the range of the data along each axis, i.e., \(\frac{\max(X_i) - \min(X_i)}{\ell}\) for \(i = 1, \dots, d\).

The density function estimator $f_{n,\ell}(\mathbf{X})$ is generalized to the multivariate case as follows:
\[
\hat f_{n,\ell}(\mathbf{x}) =
\frac{n_k}{n} \cdot \prod_{i=1}^d \frac{2m_k}{n_k (x_{i,(a_i + m_k)}^k - x_{i,(a_i - m_k)}^k)}
\]
where $\xi_{a_j}^k < x_{j} \leq \xi_{a_j + 1}^k, j = 1, \dots, d$ and $k = 1, \dots, \ell^d$ indexes the partitions, \(n_k\) is the number of samples in partition \(k\), \(m_k\) is the spacing parameter, \(a_i\) is the index of the order statistic corresponding to the interval that $X_i$ falls into and $x_{i,(r)}^k$ denote the $r$th marginal order statistic for the $i$th dimension within partition $k$. The sub-grid points for each dimension is defined as \(\xi_{i,a_i}^k = \frac{x_{i,(a_i - m_k)}^k + \dots + x_{i,(a_i + m_k - 1)}^k}{2m_k}\). Analogous to the bivariate case, the sub-grid is constructed on the empirical support of the data within each partition. Therefore, the outer boundaries for each dimension $j$ are set to the minimum and maximum marginal order statistics within that partition: $\xi_{j,0}^k = x_{j,(1)}^k \; \text{and} \; \xi_{j,n_k+1}^k = x_{j,(n_k)}^k.$

\paragraph{Remark.}
The local spacing-based density is defined only for partitions with $n_k \ge 2$,
since valid $m$-spacings require at least two sample points in the cell.
Partitions with $n_k = 0$ or $n_k = 1$ are skipped and do not contribute to 
$\hat f_{n,\ell}(X_v)$ or to the entropy estimator.

\begin{proposition}
The partitioned multivariate sample-spacing density estimator $\hat{f}_{n,\ell}(\mathbf{x})$ satisfies
\[
\lim_{n\rightarrow\infty}\int_{\mathbb{R}^d} \hat{f}_{n,\ell}(\mathbf{x}) \,d\mathbf{x} = 1.
\]
\end{proposition}

See Appendix H for the proof sketch. Now that we have established that $\hat{f}_{n,\ell}(\mathbf{x})$ is an asymptotically valid probability density function, we proceed to the main result for the estimator: proving that it is also strongly consistent. This means showing that the estimator converges almost surely to the true underlying density $f(\mathbf{x})$ as the sample size increases. The following theorem formalizes this crucial property.

\begin{theorem}[Consistency of the Partitioned Multivariate Density Estimator]\label{the_consistency_multi_density}
The partitioned multivariate density estimator is defined as:
\[
 \hat{f}_{n,\ell}(\mathbf{x}) = \frac{n_k}{n} \cdot \prod_{j=1}^{d} \widehat f_{X_j;k}(x_j),
\]
where $\widehat f_{X_j;k}(x_j)=\frac{2m_k}{n_k (x_{j,(a_j + m_k)}^k - x_{j,(a_j - m_k)}^k)}$. 
Assume that as $n \to \infty$, the number of partitions satisfies $\ell(n) \to \infty$ with $\ell(n)^d = o(n)$. Further, assume the spacing parameter $m_k$ satisfies the conditions $m_k \to \infty$ and $m_k/n_k \to 0$, which are met by the globally adopted rate of $m_k = \lfloor\sqrt{n_k} + 1/2\rfloor$. Given these conditions, $\hat{f}_{n,\ell}(\mathbf{x})$ converges almost surely and in $L^1$ to the true joint density $f(\mathbf{x})$:
\[
 \hat{f}_{n,\ell}(\mathbf{x}) \xrightarrow{a.s.} f(\mathbf{x})\quad \text{and} \quad \hat{f}_{n,\ell}(\mathbf{x}) \xrightarrow{L^1} f(\mathbf{x}).
\]
\end{theorem}
See Appendix I for the proof sketch.
\subsection{Multivariate Joint Entropy Estimator}

The multivariate joint entropy estimator is defined as the plug-in estimator based on the partitioned multivariate density:
\[
    \hat H_{n,\ell}(\mathbf{X}) = -\frac{1}{n} \sum_{v=1}^{n} \log \hat{f}_{n,\ell}(\mathbf{X}_v).
\]
For a given data point $\mathbf{X}_v$ that falls into a sub-grid cell indexed by $\mathbf{a} = (a_1, \dots, a_d)$ within partition $k$, the log-density term is expanded as:
\[
-\log \hat{f}_{n,\ell}(\mathbf{X}_v) = -\log\left(\frac{n_k}{n} \cdot \prod_{j=1}^{d} \frac{2m_k}{n_k \Delta x_{j, a_j}^k} \right),
\]
where $\Delta x_{j, a_j}^k = x_{j,(a_j + m_k)}^k - x_{j,(a_j - m_k)}^k$. The final estimator is the average of these terms over all $n$ data points.

Building on Theorem~\ref{the_consistency_multi_density}, the convergence of the entropy estimator can now be established. The proof follows analogously to that of the bivariate entropy consistency. This final theoretical result is stated in the following theorem.

\begin{theorem}[Consistency of the Entropy Estimator $\hat{H}_{n,\ell}$]\label{thm_consistency_entropy}
Assume that the true density $f$ satisfies $\log f \in L^{1+\delta}$ for some $\delta \in (0,1)$, $m_k = \lfloor \sqrt{n_k} + 1/2 \rfloor$, with $n_k > 1$, and that the number of partitions satisfies $\ell(n)^d = o(n)$. Then, the entropy estimator $\hat{H}_{n,\ell}$
converges in $L^1$ to the true entropy $H(f) := \mathbb{E}[-\log f(X)]$ as $n \to \infty$:
\[
\hat{H}_{n,\ell} \xrightarrow{L^1} H(f), \quad \text{and} \quad \hat{H}_{n,\ell} \xrightarrow{\text{a.s.}} H(f).
\]
\end{theorem}

See Appendix J for the proof sketch. Building upon the $L^1$ and almost sure convergence of the PSS entropy estimator 
established in Theorem~\ref{thm_consistency_entropy}, we can immediately obtain the consistency of the 
corresponding mutual information estimator.

\begin{corollary}[Consistency of the PSS Mutual Information Estimator]
Let $X \in \mathbb{R}^{d_X}$ and $Y \in \mathbb{R}^{d_Y}$ be random vectors. 
Define the PSS mutual information estimator as
\[
\hat I_{n,\ell}(X;Y) \;=\; \hat H_{n,\ell}(X) + \hat H_{n,\ell}(Y) - \hat H_{n,\ell}(X,Y).
\]
Under the same assumptions as Theorem~\ref{thm_consistency_entropy}, 
\[
\hat I_{n,\ell}(X;Y) \xrightarrow{L^1} I(X;Y),
\quad
\hat I_{n,\ell}(X;Y)\xrightarrow{a.s.}I(X;Y).
\]
\end{corollary}

\subsection{Data-Driven Parameter Selection via Cross-Validation}
\label{sec:cv_selection}

A critical challenge in nonparametric estimation is the selection of the smoothing parameter---in our case, the number of partitions per dimension, $\ell$. While theoretical analyses for grid-based estimators suggest that the partition size should scale with the sample size, the optimal rate typically depends on unknown properties of the true density $f$ (e.g., roughness or derivatives). To address this in a practical, data-driven manner without relying on unverifiable assumptions, we propose a Likelihood-based Cross-Validation strategy.

Let the dataset $\mathcal{D} = \{\mathbf{x}_1, \dots, \mathbf{x}_n\}$ be partitioned into $K$ disjoint folds. For each fold $k \in \{1, \dots, K\}$, let $\mathcal{V}_k$ denote the validation set and $\mathcal{T}_k = \mathcal{D} \setminus \mathcal{V}_k$ the training set. We select the optimal parameter $\ell^*$ by minimizing the average negative log-likelihood on the hold-out data:
\begin{equation}
    \label{eq:cv_criterion}
    \ell^* = \operatorname*{argmin}_{\ell } \left\{ - \frac{1}{K} \sum_{k=1}^{K} \frac{1}{|\mathcal{V}_k|} \sum_{\mathbf{x} \in \mathcal{V}_k} \log \hat{f}_{n, \ell}^{(\mathcal{T}_k)}(\mathbf{x}) \right\},
\end{equation}
where $\hat{f}_{n, \ell}^{(\mathcal{T}_k)}$ represents the PSS density estimator constructed using only the training subset $\mathcal{T}_k$.

\paragraph{Connection to the Entropy Objective.}
The proposed cross-validation criterion is not merely a heuristic but is intrinsically linked to the objective of entropy estimation. Recall that the plug-in entropy estimator is defined as the empirical average of the negative log-density:
\begin{equation}
    \hat{H}_{n, \ell} = - \frac{1}{n} \sum_{i=1}^{n} \log \hat{f}_{n, \ell}(\mathbf{x}_i) \, .
\end{equation}
Comparing Eq.~\eqref{eq:cv_criterion} with the definition of $\hat{H}_{n, \ell}$, it is evident that the validation loss $\mathcal{L}_{CV}$ takes an identical form to the entropy estimator itself. 
By the law of large numbers, the expectation of this loss converges to the cross-entropy between the true density $f$ and the estimator $\hat{f}$:
\begin{equation}
    \mathbb{E}_{\mathbf{x} \sim f} \left[ - \log \hat{f}_{n, \ell}(\mathbf{x}) \right] = H(f) + D_{KL}(f \,||\, \hat{f}_{n, \ell}) \, .
\end{equation}
Since the true entropy $H(f)$ is constant with respect to $\ell$, minimizing the CV loss is asymptotically equivalent to minimizing the Kullback-Leibler (KL) divergence $D_{KL}(f \,||\, \hat{f}_{n, \ell})$. This structural alignment ensures that the hyperparameter $\ell$ is optimized specifically for the goal of accurate information estimation, providing a coherent and principled framework for automated model selection.

\subsection{Practical Implementation Note}

The asymptotic condition $\ell(n)^d = o(n)$ in our theoretical analysis is a sufficient but not a necessary requirement for consistency. In practical applications, the PSS estimator evaluates $\hat f_{n,\ell}(X_v)$ only at observed sample points; partitions with $n_k = 0$ or with insufficient spacing do not contribute to the entropy computation. As a result, the estimator is governed by the number of \emph{occupied} partitions rather than by the full $\ell^d$ grid.

Real datasets, including the whitened 14-dimensional EEG signals used in our experiments, typically concentrate on a lower-dimensional manifold. Only a small fraction of the $\ell^d$ hyperrectangles actually contain data, and the entropy estimator acts only on these populated regions. Consequently, practical choices such as $\ell^\ast = 12$ remain stable and do not conflict with the asymptotic condition, since the estimator adapts to the effective geometry and occupancy pattern of the data.

\vspace{0.5em}

\begin{algorithm}[H]
\caption{PSS Joint Entropy Estimation}
\label{alg:pss_entropy}
\begin{algorithmic}[1]
    \Require Samples $S=\{\mathbf{x}_v\}_{v=1}^n \subset \mathbb{R}^d$, partition parameter $\ell$.
    \Ensure Estimated joint entropy $\hat H_{n,\ell}$.
    
    \State Partition each axis into $\ell$ equal-width intervals to form $\ell^d$ cells $\{P_k\}$.
    \State $\mathcal{S} \gets 0$ \Comment{Initialize log-density sum}
    
    \Statex \Comment{\textit{Pre-computation step}}
    \For{each occupied partition $P_k$ with $n_k \ge 2$}
        \State Sort points in $P_k$ marginally.
        \State Compute spacings $\Delta x^k_{j,a_j} = x^k_{j,(a_j+m_k)} - x^k_{j,(a_j-m_k)}$ using $m_k=\lfloor\sqrt{n_k}\rfloor$.
    \EndFor
    
    \Statex \Comment{\textit{Entropy estimation step}}
    \For{$v=1$ to $n$}
        \State Identify partition $P_k$ containing $\mathbf{x}_v$.
        \If{$n_k < 2$ or any $\Delta x^k_{j,a_j}=0$} 
            \State \textbf{continue} \Comment{Skip points with undefined density}
        \EndIf
        
        \State Calculate local log-density contribution:
        \[
            L_v = \log\left(\frac{n_k}{n}\right) + \sum_{j=1}^d \log\left(\frac{2 m_k}{n_k \,\Delta x^k_{j,a_j}}\right)
        \]
        \State $\mathcal{S} \gets \mathcal{S} + L_v$
    \EndFor

    \State $\hat H_{n,\ell} = -\frac{1}{n} \mathcal{S}$
    \State \Return $\hat H_{n,\ell}$
\end{algorithmic}
\end{algorithm}

\newpage
\section{Numerical Experiments}
\label{sec:experiments}

To comprehensively evaluate the performance and practical utility of the proposed Partitioned Sample-Spacing (PSS) estimator, we conducted extensive experiments on both synthetic and real-world datasets. Our experimental design focuses on validating three key properties: (1) robustness to skewed distributions, (2) computational scalability in large-sample regimes, and (3) reliability in high-dimensional downstream tasks.

\paragraph{Hyperparameter Selection.}
For hyperparameter selection---specifically, the partition parameter $\ell$ for PSS and the number of neighbors $k$ for kNN---we adopted distinct protocols depending on the availability of ground truth:
\begin{itemize}
    \item \textbf{Synthetic Benchmarks (Oracle Tuning):} In scenarios where the true probability density and entropy $H(f)$ are known, we selected hyperparameters that minimized the MSE. This allows us to assess the empirical bound of the estimator's accuracy and convergence properties, independent of hyperparameter selection error.
    \item \textbf{Real-World Applications (Data-Driven Tuning):} For real-world tasks such as ICA and feature selection, where the ground truth is inaccessible, we employed the Likelihood-based Cross-Validation(CV) strategy proposed in Section~4.3. This evaluates the estimators in a realistic setting, demonstrating their stability and reliability when fully automated, data-driven tuning is required.
\end{itemize}

\subsection{Synthetic Data Experiments}

We compare PSS against $k$NN estimators (KL, KSG), copula based CADEE \citep{copula}, and uniformizing-mapping versions (UM-tKL, UM-tKSG) \citet{AoLi2022}. For both CADEE and the UM-based estimators, we adapted the authors' official implementations to fit our experimental framework. We do not include general NF-based entropy estimators as baselines. As shown in \citep{AoLi2022}, pure normalizing-flow estimators are
typically less accurate and less stable than the normalized truncated $k$NN estimators. To ensure a fair comparison, the key hyperparameter for each method---the partition parameter $\ell$ for PSS and the number of neighbors $k$ for kNN estimators---was chosen to minimize the empirical Mean Squared Error (MSE) in each setting. All results are averaged over 100 trials and reported as rooted mean squared error (RMSE) relative to the true entropy. 

\paragraph{Computational Complexity Analysis.}
The standard $k$NN estimator has a naive time complexity of $O(N^{2}d)$ in $d$ dimensions, since all pairwise distances are computed. 
Tree-based methods (k-d tree, Ball tree) reduce query cost to $O(Nd\log N)$ after an $O(Nd\log N)$ build phase, but their efficiency deteriorates with dimension and skewness. 
Normalizing flow (NF) estimators \citep{AoLi2022} require training a deep generative model with $O(Nd)$ forward/backward passes per epoch, repeated over many epochs, leading to substantial computational overhead. 
In contrast, PSS assigns $N$ samples to $\ell^d$ partitions in $O(Nd)$, sorts within partitions in about $O(dN \log(N/\ell^d))$, and computes entropy in $O(Nd)$. 
As $\ell$ grows with $N$, the factor $N/\ell^d$ shrinks and the cost approaches $O(Nd)$, yielding favorable scalability. 
CADEE \citep{copula} runs in reasonable time for small problems but its runtime grows rapidly with $N$ and $d$ relative to PSS, and its accuracy degrades in higher dimensions.

\newpage
\begin{figure}[h]
    \centering
    \subfloat[]{
        \includegraphics[width=0.23\textwidth]{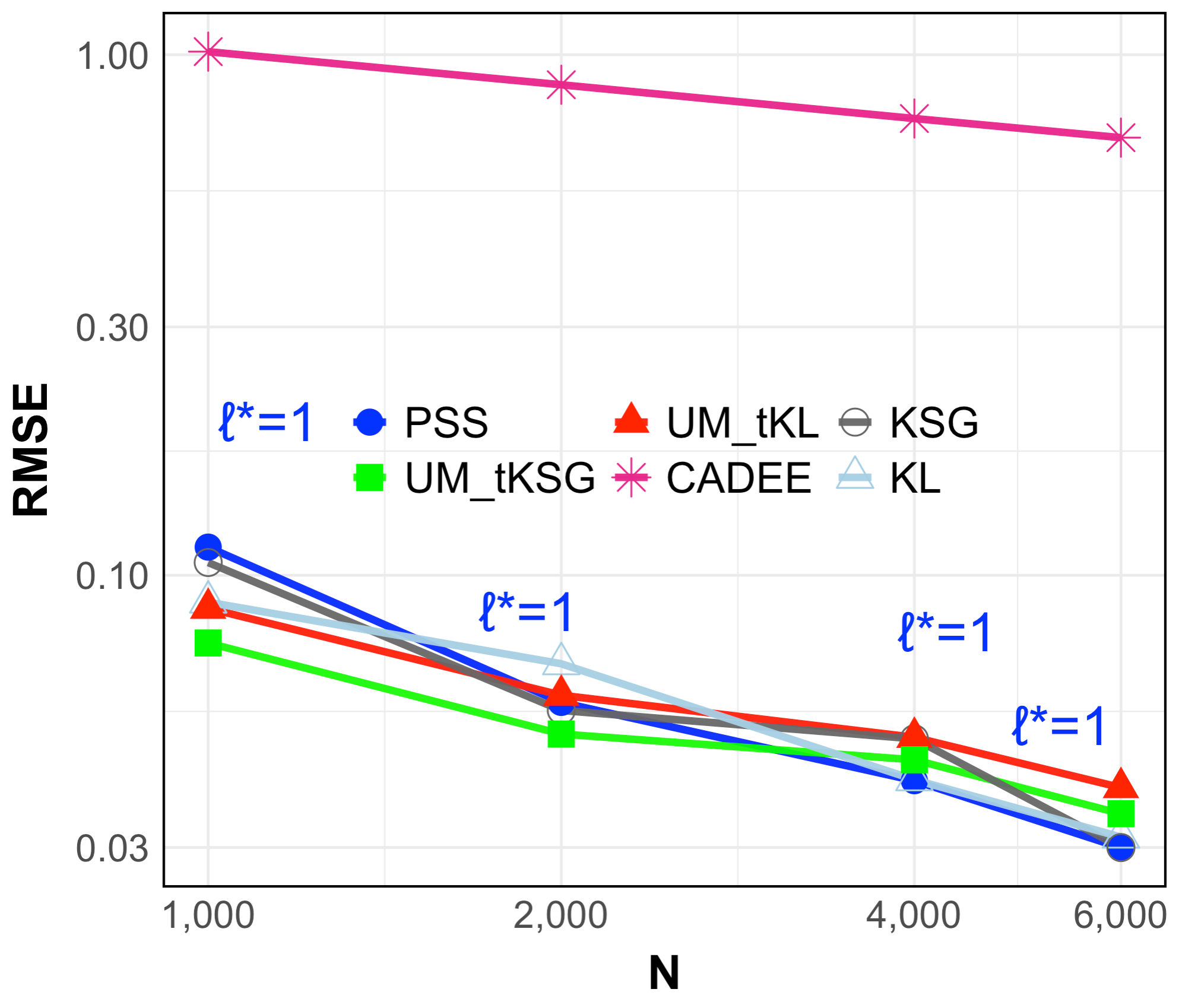}
    }
    \hfill
    \subfloat[]{
        \includegraphics[width=0.23\textwidth]{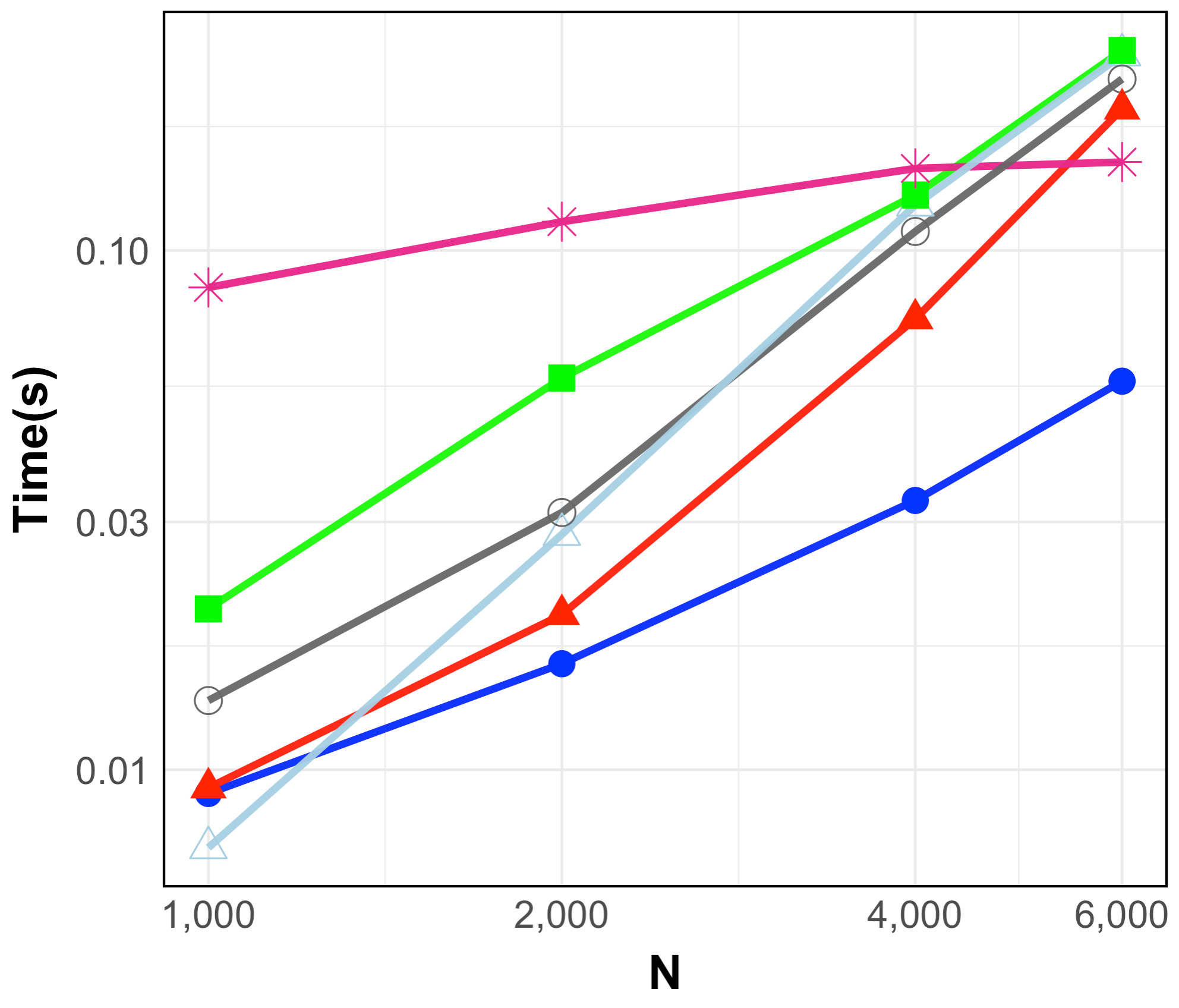}
    }
    \hfill
    \subfloat[]{
        \includegraphics[width=0.23\textwidth]{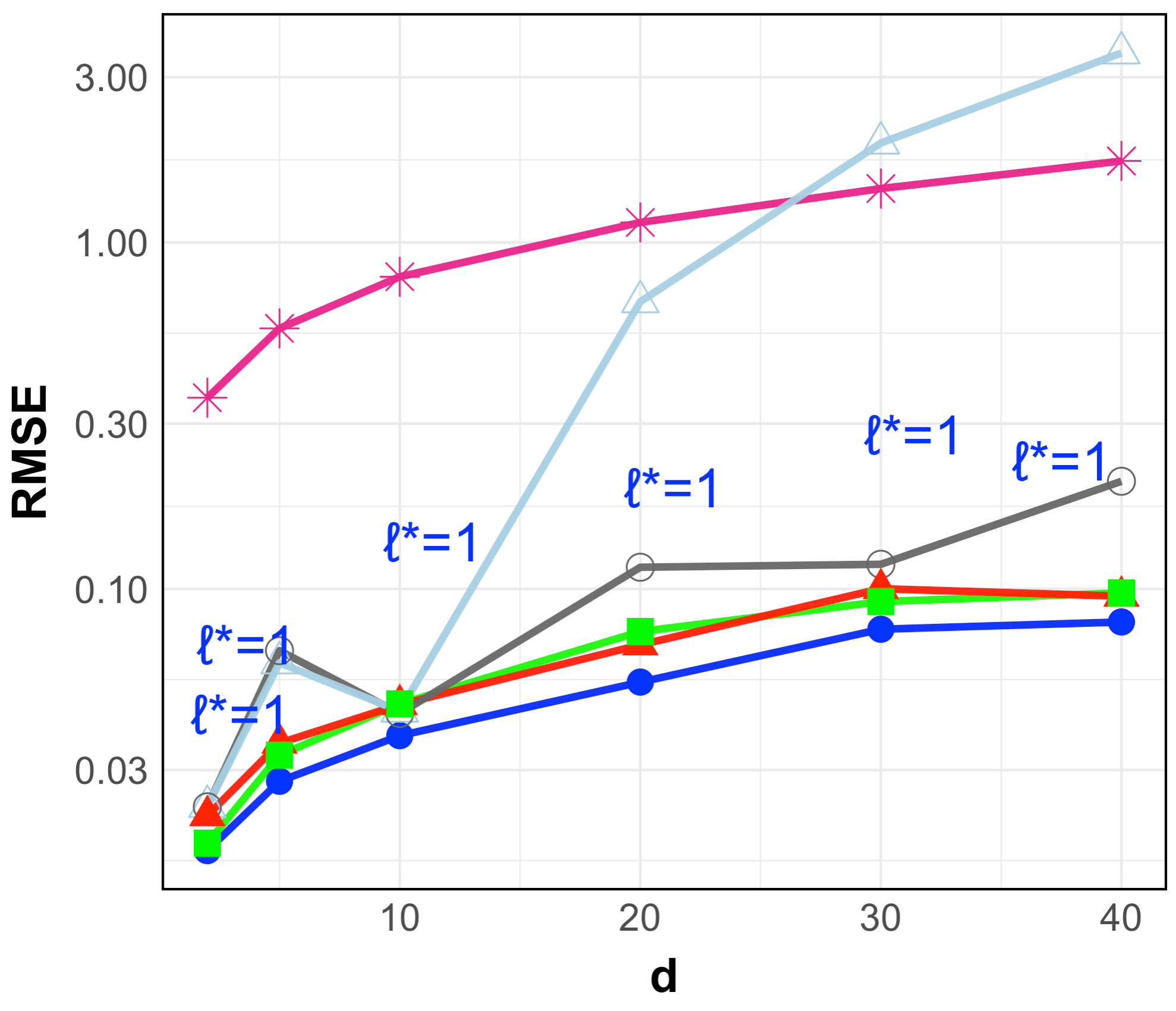}
    }
    \hfill
    \subfloat[]{
        \includegraphics[width=0.23\textwidth]{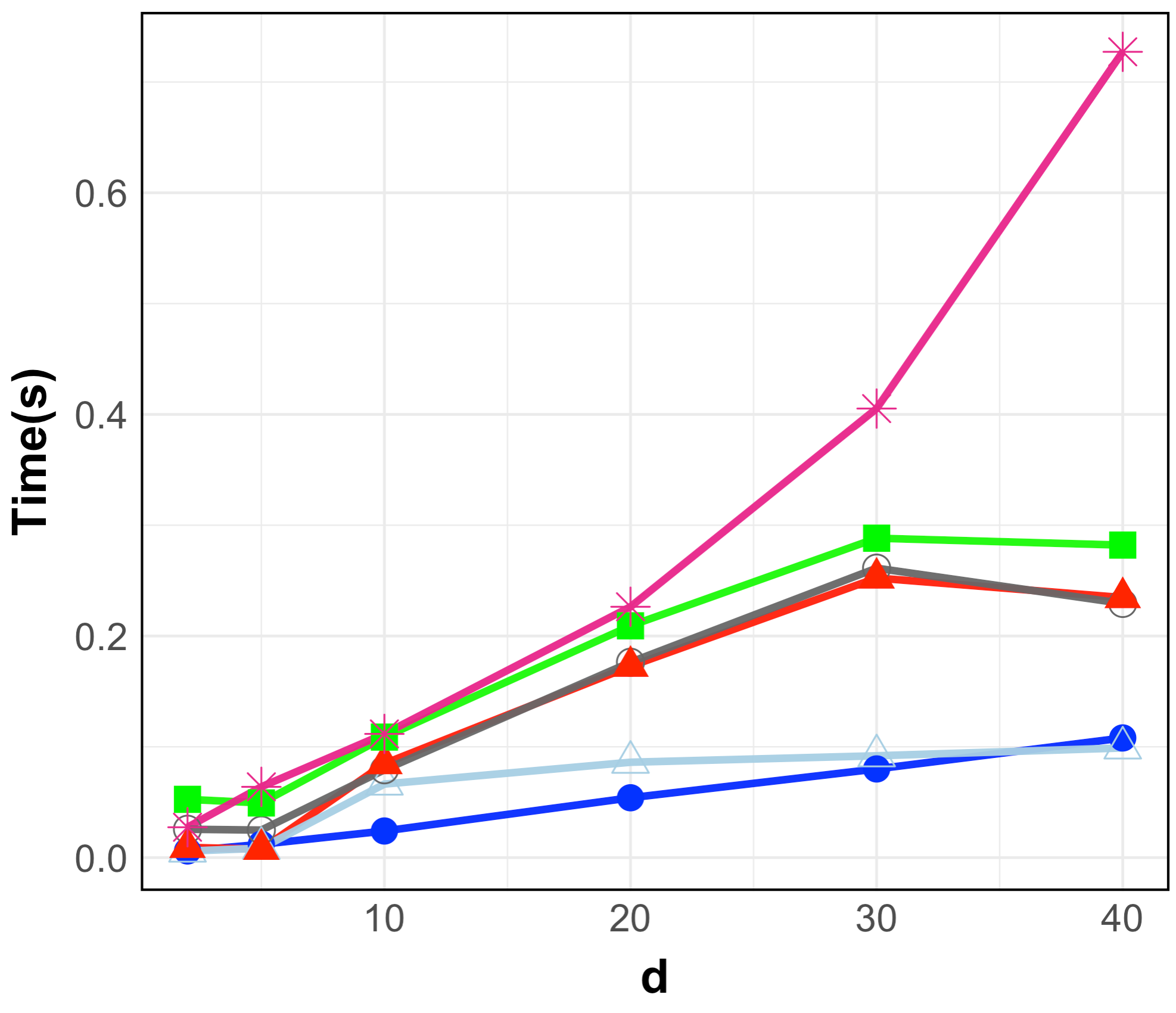}
    }
    \caption{RMSE (a, c) and runtime (b, d) for entropy estimators under the Normal distribution 
with $\rho=0$. Panels (a, b): $d=10$, varying sample size $N$. 
Panels (c, d): $N=3000$, varying dimension.}
    \label{fig:rmse}
\end{figure}

\begin{figure}[h]
    \centering
    \subfloat[]{
        \includegraphics[width=0.23\textwidth]{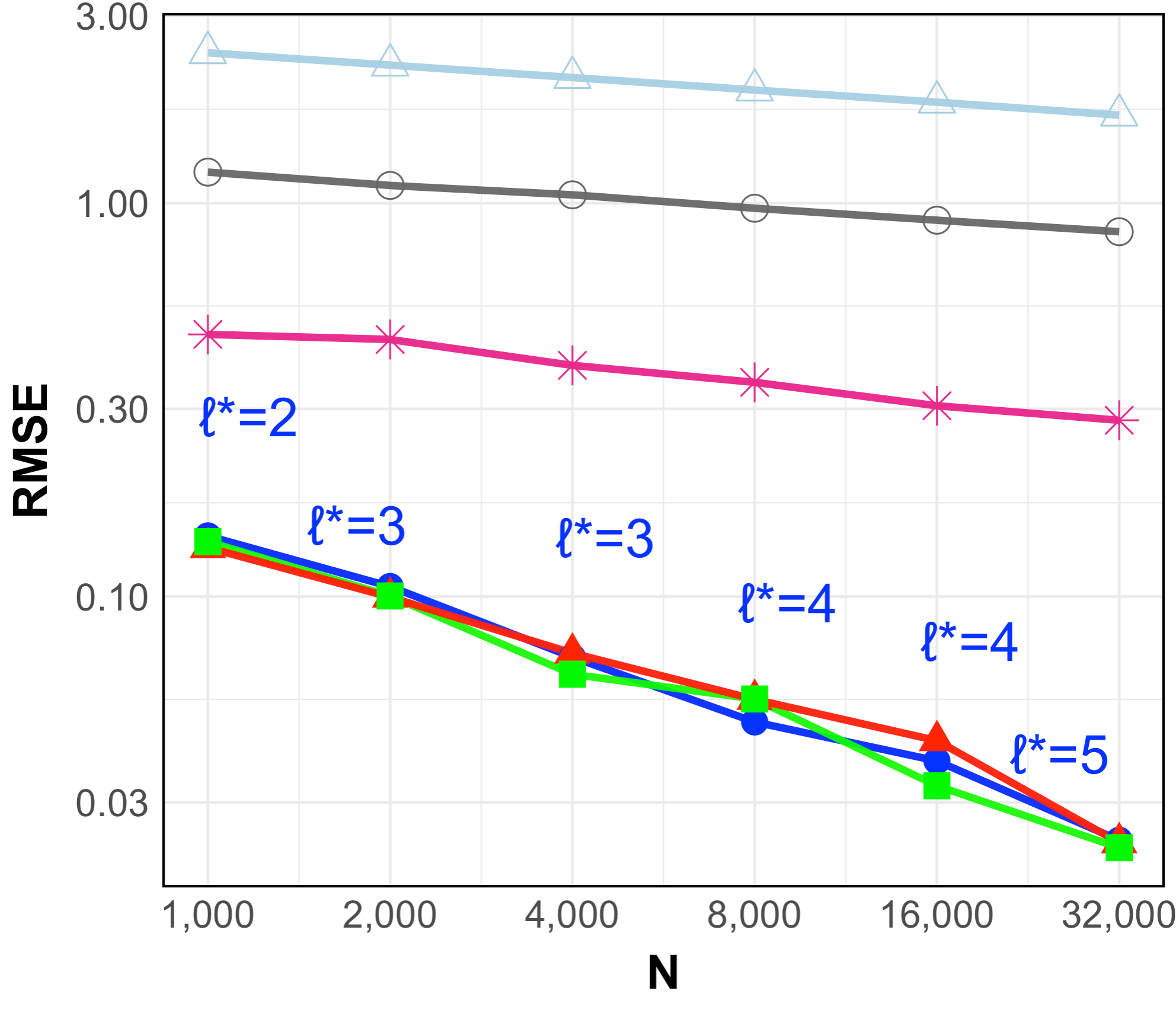}
    }
    \hfill
    \subfloat[]{
        \includegraphics[width=0.23\textwidth]{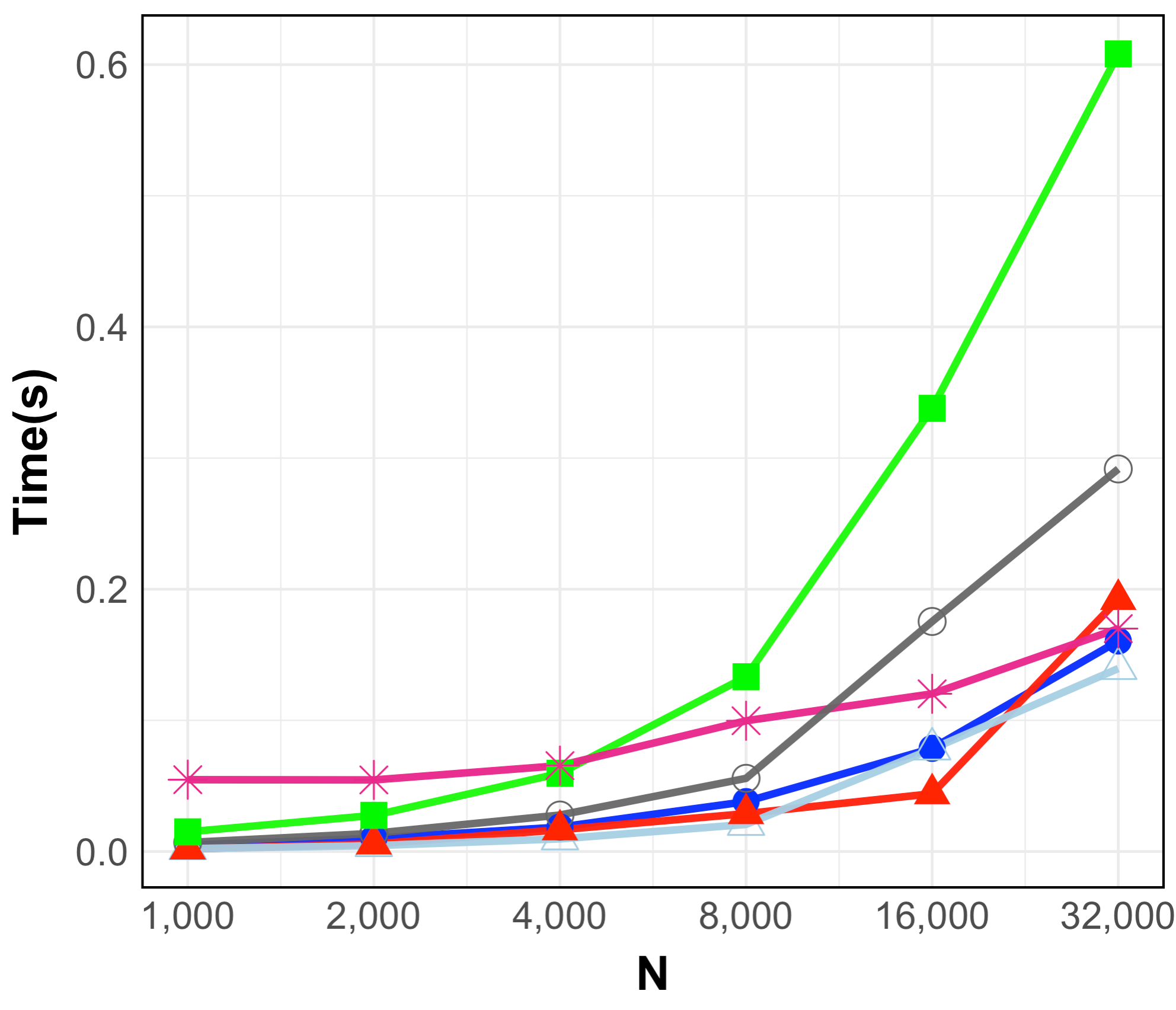}
    }
    \hfill
    \subfloat[]{
        \includegraphics[width=0.23\textwidth]{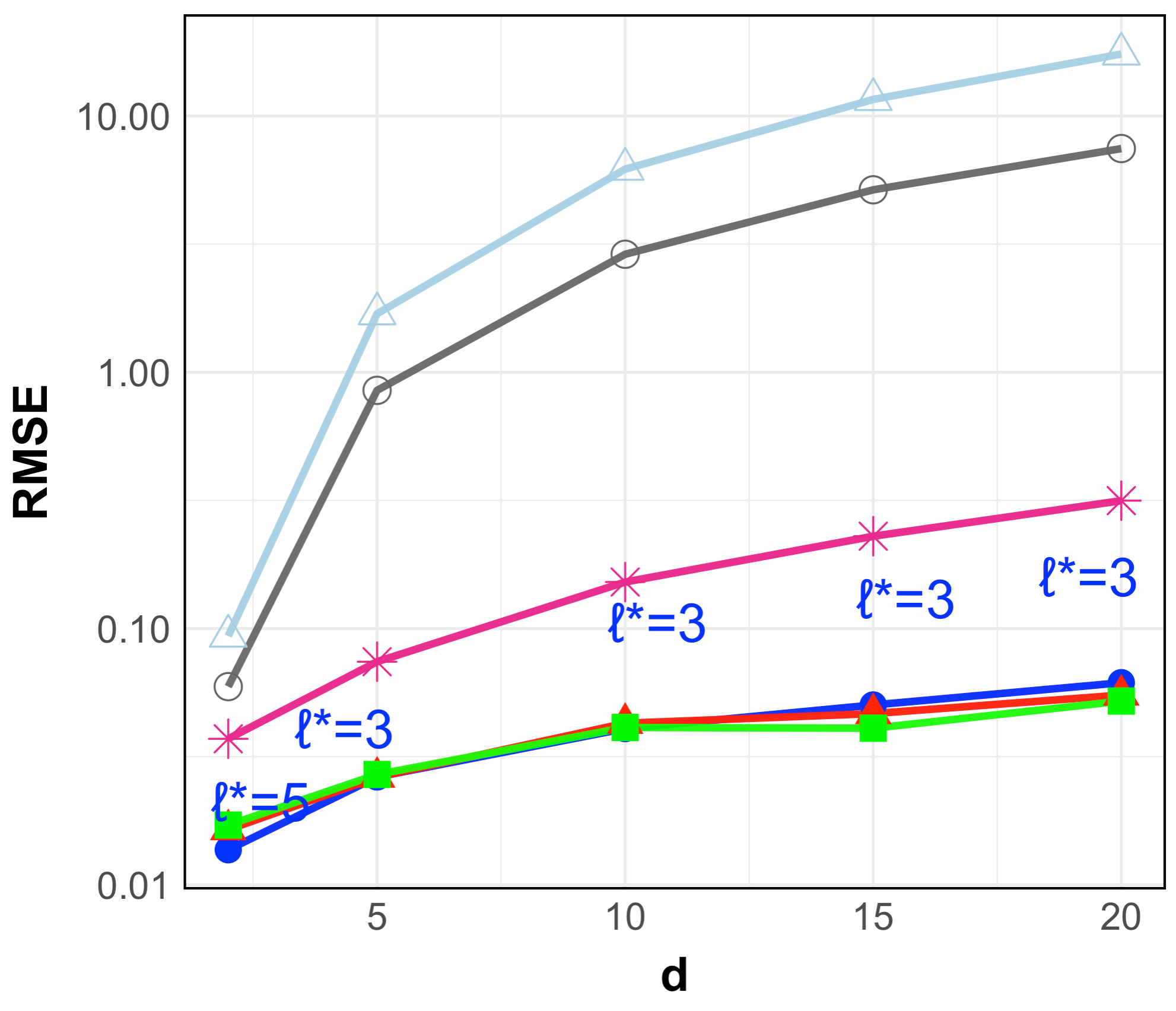}
    }
    \hfill
    \subfloat[]{
        \includegraphics[width=0.23\textwidth]{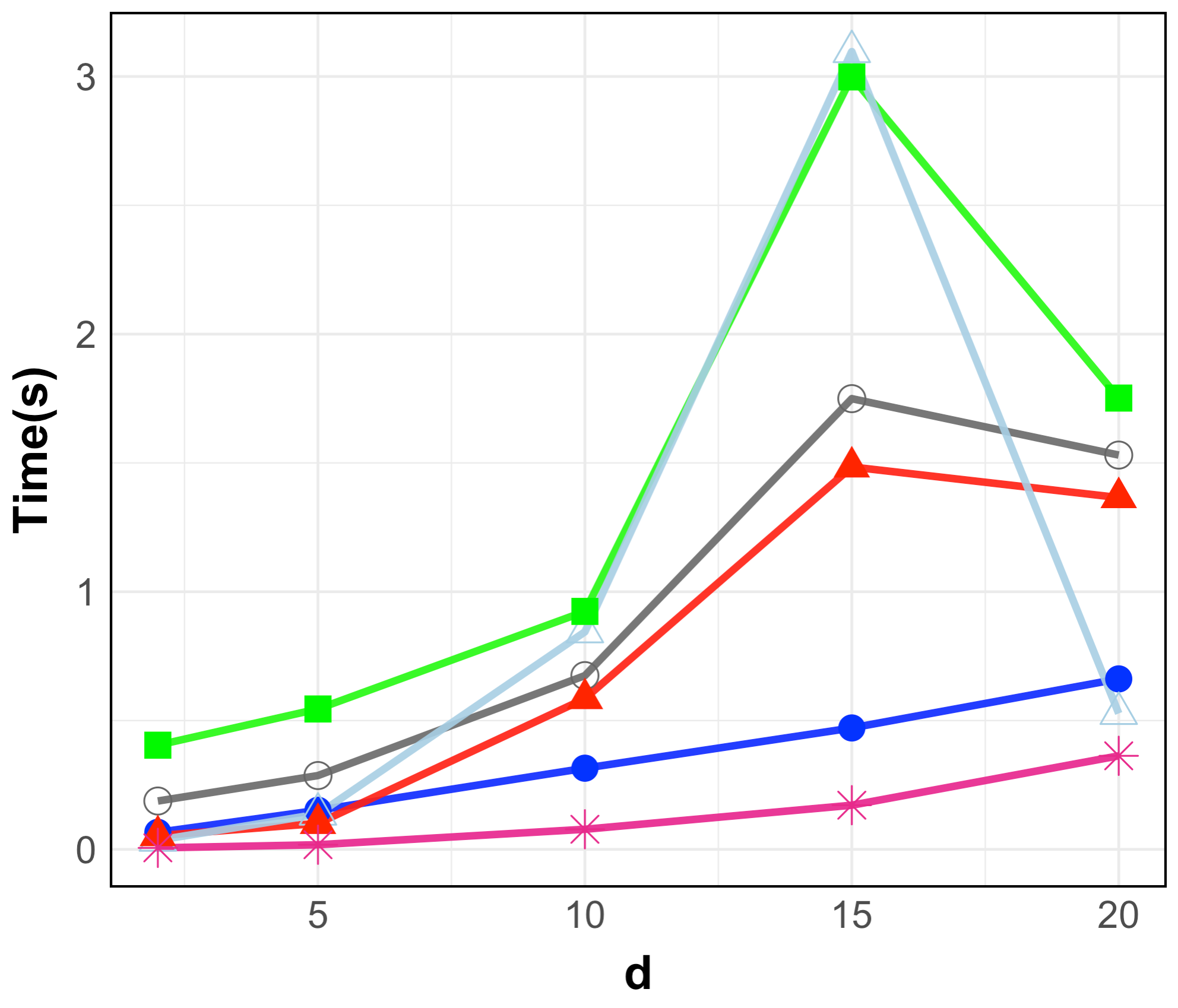}
    }
    \caption{RMSE (a, c) and runtime (b, d) for entropy estimators under the Gamma distribution with shape = 0.4, scale = 0.3. Panels (a, b): $d= 5$, varying sample size $N$. Panels (c, d): $\rho=0$, $N=30,000$, varying dimension.}
\label{fig:rmse2}
\end{figure}

\begin{figure}[h]
    \centering
    \subfloat[]{
        \includegraphics[width=0.23\textwidth]{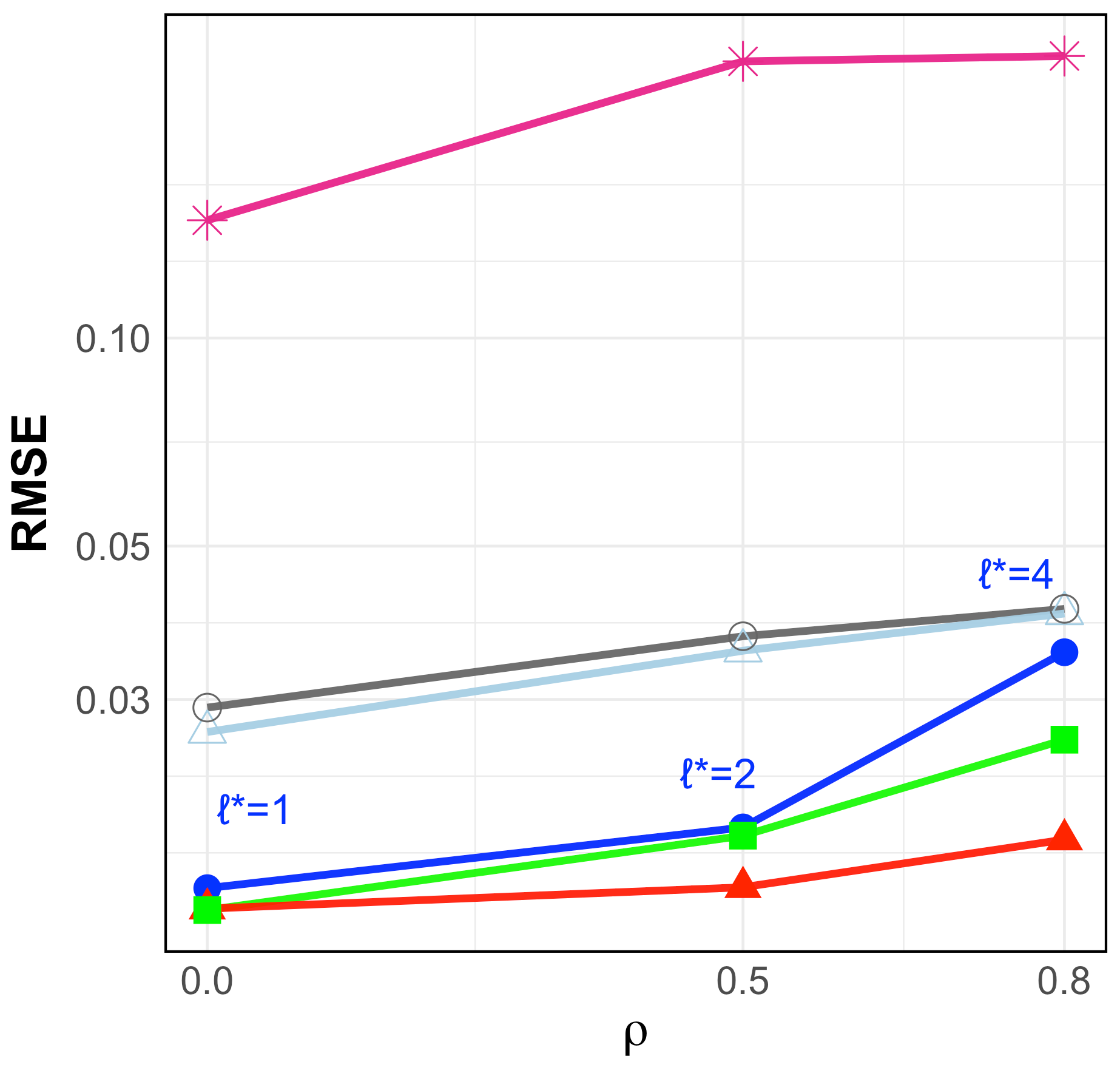}
    }
    \hfill
    \subfloat[]{
        \includegraphics[width=0.23\textwidth]{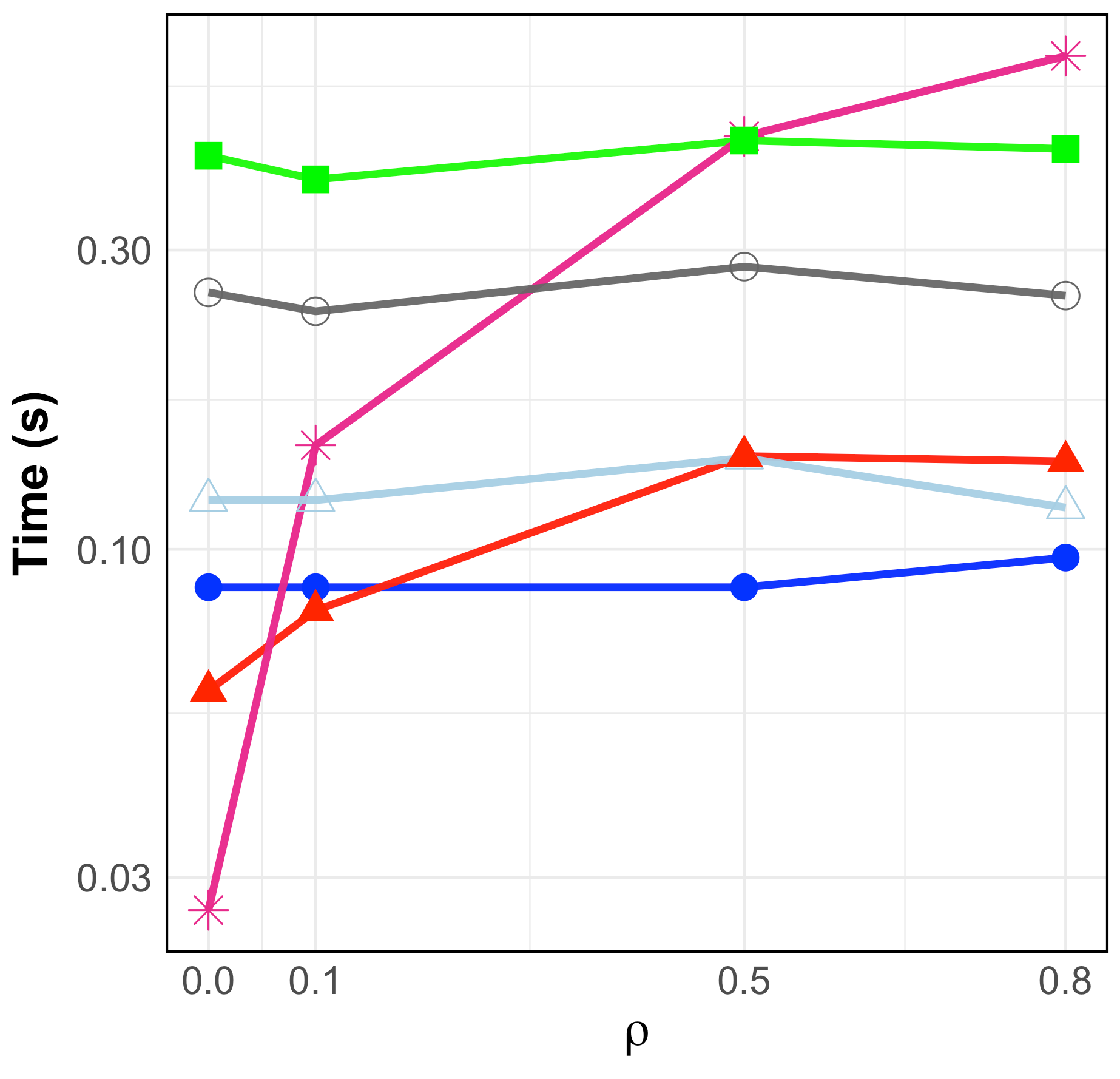}
    }
    \hfill
    \subfloat[]{
        \includegraphics[width=0.23\textwidth]{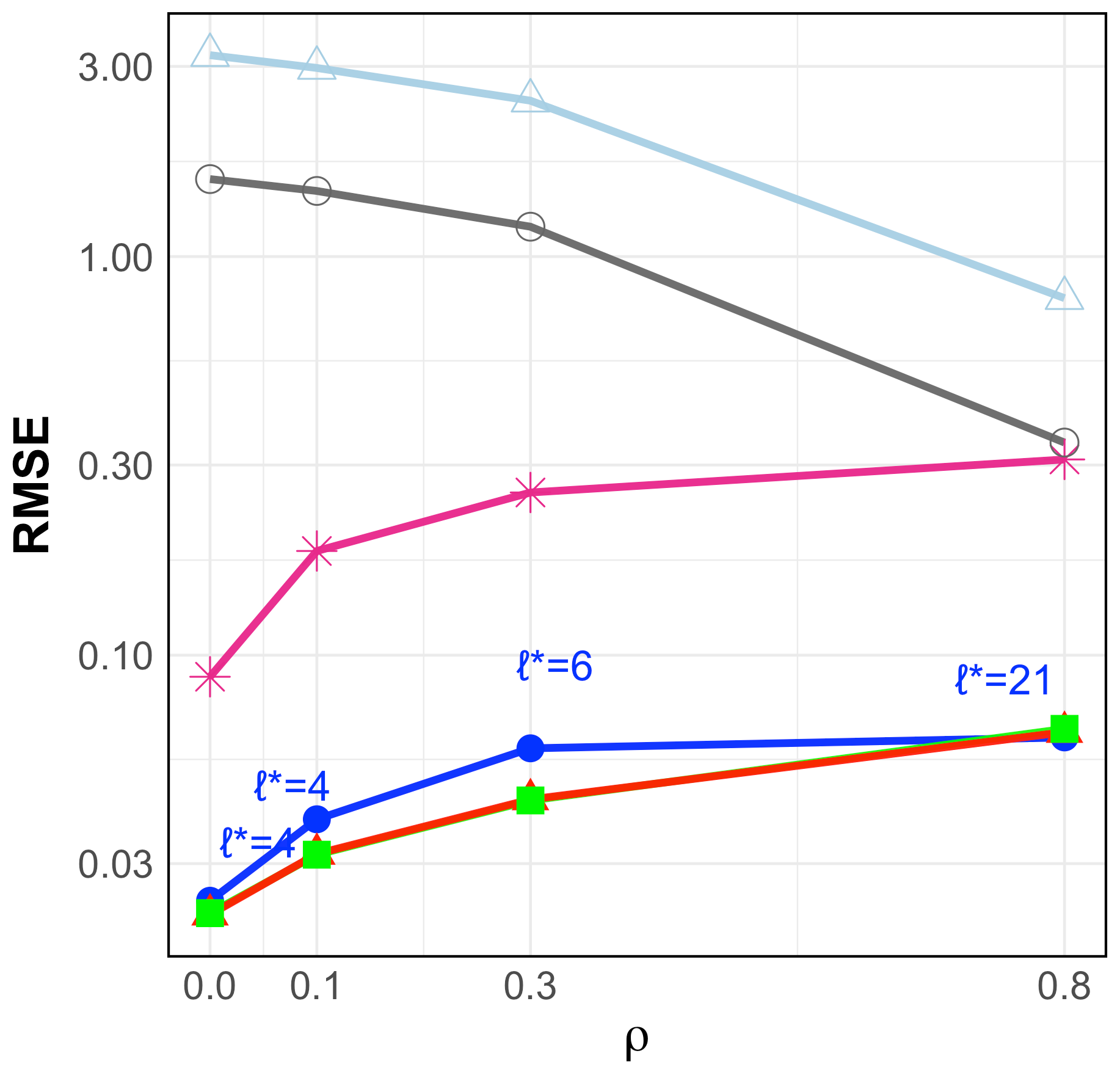}
    }
    \hfill
    \subfloat[]{
        \includegraphics[width=0.23\textwidth]{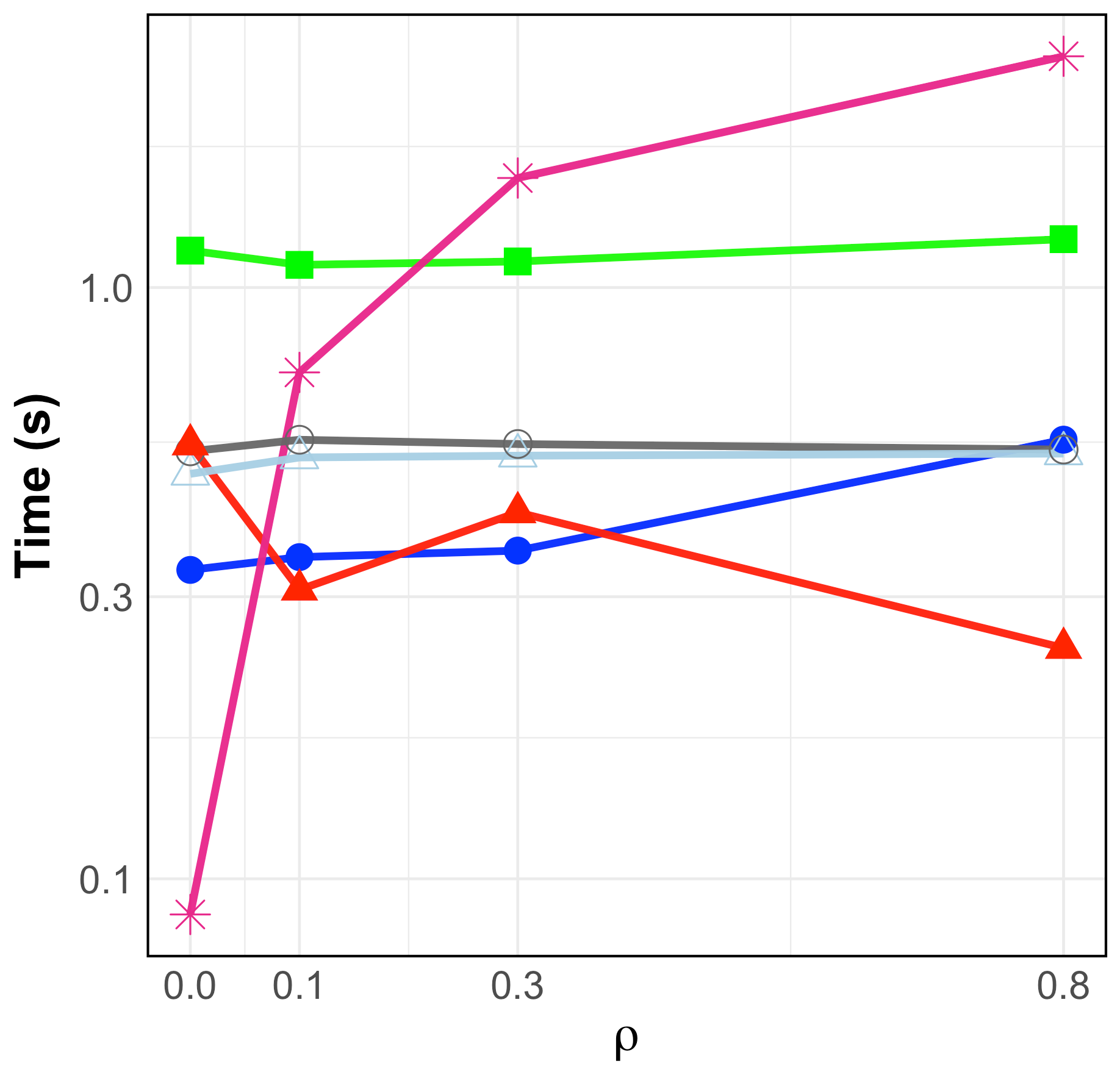}
    }
    \caption{RMSE (a, c) and runtime (b, d) for entropy estimators with varying correlation coefficient $\rho$. 
Panels (a, b): Normal distribution ($d=5, N=20,000$). 
Panels (c, d): Gamma distribution (shape $=0.4$, scale $=0.3$) with $d=7, N=50,000$.}
\label{fig:rmse_rho}
\end{figure}

\paragraph{Multivariate Normal.}
We first consider multivariate normal data with $\rho=0$. 
Figure~\ref{fig:rmse} reports RMSE and runtime across varying sample sizes ($d=10$) and dimensions ($N=3000$). 
In terms of accuracy (Figure~\ref{fig:rmse}a, c), PSS consistently outperforms traditional baselines (KL, KSG) and CADEE, with the performance gap widening as $d$ grows. 
Notably, PSS achieves accuracy competitive with state-of-the-art normalizing flow-based estimators (UM-tKL, UM-tKSG), effectively matching their performance even in higher dimensions ($d=40$) without the need for auxiliary density modeling.
In terms of runtime (Figure~\ref{fig:rmse}b, d), PSS scales nearly linearly with $N$ and remains substantially faster than CADEE and the UM-variants, indicating that PSS alleviates several dimensionality-related challenges in Gaussian settings while maintaining strong computational efficiency.

\paragraph{Multivariate Gamma.}
We next evaluate Gamma marginals combined through a Gaussian copula with $\rho=0$, shape $0.4$, and scale $0.3$ (Figure~\ref{fig:rmse2}). 
In all Gaussian–copula experiments, the oracle joint entropy is computed analytically as
\[
H(X) \;=\; \sum_{j=1}^d h(X_j) + \tfrac{1}{2}\log\det R,
\]
where $h(X_j)$ is the marginal Gamma entropy and $R$ is the copula correlation matrix. 
In terms of accuracy, PSS exhibits remarkable robustness to skewed distributions and semi bounded supports. 
It significantly outperforms traditional baselines (KL, KSG) and CADEE, which suffer from severe boundary bias in this setting. 
Crucially, PSS matches or exceeds the performance of normalizing flow-based methods (UM-tKL, UM-tKSG), demonstrating that the partitioning strategy effectively captures local density structures without relying on learned transformations.
Runtime comparisons further highlight the advantage of PSS; it scales linearly with $N$ and remains substantially faster than CADEE, while avoiding the computational overhead associated with training flow-based models.

\paragraph{Robustness under Correlation.}
We further investigate the estimator's sensitivity to dependence strength by varying the correlation coefficient $\rho \in [0, 0.8]$ (Figure~\ref{fig:rmse_rho}). 
In terms of accuracy (Figure~\ref{fig:rmse_rho}a, c), PSS exhibits remarkable stability, maintaining low errors even as traditional baselines degrade rapidly with increasing correlation.
Notably, in the high-correlation Gamma setting ($\rho=0.8$), PSS surpasses even the normalizing flow-based variants without requiring learned transformations.
This robustness is directly attributable to the adaptive behavior of the optimal partition parameter $\ell^*$.
We observe that a larger $\ell^*$ is required to minimize error as correlation strengthens (e.g., increasing from 4 to 21 in the Gamma case), indicating that a finer partitioning is essential to capture the probability mass concentrating in narrow regions.
Regarding computational cost (Figure~\ref{fig:rmse_rho}b, d), PSS displays a flat runtime profile, confirming that its complexity remains invariant to the dependence structure, whereas CADEE suffers significant slowdowns.
However, we note a structural limitation as $\rho \to 1$; while PSS handles $\rho=0.8$ robustly, the transition to a fully degenerate distribution on a lower-dimensional manifold remains a fundamental challenge for full-dimensional density estimators.

\subsection{Real-World Data Experiments}
\subsubsection{ICA Evaluation Experiment}
To demonstrate the robustness and reliability of our proposed Partitioned Sample-Spacing (PSS) estimator on real-world, moderate-dimensional data, we conducted an experiment evaluating the performance of Independent Component Analysis (ICA). The procedure is summarized as follows.

\paragraph{1. Data and Preprocessing.}
We used the UCI \texttt{EEG Eye State} dataset\citep{eeg_eye_state_264}, a 14-dimensional time-series signal comprising $N = 14{,}980$ samples and $d = 14$ channels. The dataset is a continuous recording from a single subject. To prepare the data for entropy estimation---which assumes independent and identically distributed (i.i.d.) samples---we first pre-processed the entire time series via whitening, which zero-centers the data and applies a linear transformation to decorrelate the features. Each row of the transformed data matrix was then treated as an i.i.d. sample from a stationary $d$-dimensional distribution. This simplification ignores temporal dependence and is adopted solely for evaluating entropy estimators on real-world correlated data.

\paragraph{2. Analysis Procedure (ICA).}
We applied the \texttt{FastICA} algorithm to the pre-whitened data ($\mathbf{X}_w$) to separate it into a set of statistically independent components ($\mathbf{Y}$). All computations were performed using scikit-learn 1.5.2 (\texttt{FastICA}) on an Apple M3 processor.

\paragraph{3. Evaluation Metric (Total Correlation).}
To quantify the performance of ICA, we measure the Total Correlation (TC), which represents the total statistical dependence in a random vector $\mathbf{X} = (X_1, \ldots, X_d)$. It is defined as the sum of the marginal entropies minus the joint entropy:
\begin{equation*}
    \mathrm{TC}(\mathbf{X}) = \sum_{i=1}^{d} H(X_i) - H(\mathbf{X}) \, .
\end{equation*}
TC is non-negative, with $\mathrm{TC} = 0$ indicating perfect independence. A successful ICA should yield components with a TC value close to zero. We measured the TC of the data before applying ICA ($\texttt{TC}_{\text{before}}$) and after ($\texttt{TC}_{\text{after}}$). An estimator’s reliability is judged by its ability to report a small, positive $\texttt{TC}_{\text{after}}$ value.

\paragraph{4. Hyperparameter Tuning.} 
For a fair comparison, the key hyperparameter for each estimator was selected via 3-fold cross-validation ($K=3$) on the whitened data. As proposed in the discussion section, we employed a principled, likelihood-based validation strategy for PSS: selecting the number of partitions ($\ell$) that minimized the average negative log-density on the held-out folds. This criterion is directly aligned with the entropy estimation objective, as the estimator itself is defined as an empirical average of the negative log-density. This single, globally-optimized hyperparameter, denoted $\ell^*$, was then used for all subsequent PSS-based calculations, including both the joint and the marginal entropy estimations. For KL/KSG, we minimize the held-out negative log-density estimated.

\begin{figure}[h]
\centering
\begin{tikzpicture}
\begin{axis}[
    title={\textbf{Change in Total Correlation (TC) Before and After ICA}},
    width=0.8\textwidth,
    height=8cm,
    ylabel={TC (nats)},
    xlabel={Entropy Estimator},
    symbolic x coords={PSS, KSG, tKSG, KL, tKL},
    xtick=data,
    ymin=-12, ymax=14,
    enlarge x limits=0.15,
    legend style={at={(0.5,-0.25)}, anchor=north, legend columns=-1},
    ymajorgrids=true,
    grid style=dashed,
    extra y ticks={0},
    extra y tick style={grid=major, grid style={solid, black, thick}},
]
% Before (hollow)
\addplot[only marks, mark=o, blue, mark size=3pt, thick] coordinates {
    (PSS, 8.197) (KSG, 11.741) (tKSG, 11.519) (KL, 10.004) (tKL, 7.954)
};
\addlegendentry{$\texttt{TC}_{\text{before}}$}

% After (filled, valid)
\addplot[only marks, mark=*, blue, mark size=3pt, thick] coordinates {
    (PSS, 1.439) (KSG, 8.499) (tKSG, 8.084)
};
\addlegendentry{$\texttt{TC}_{\text{after}}$ (Nonnegative)}

% After (filled, invalid)
\addplot[only marks, mark=*, red, mark size=3pt, thick] coordinates {
    (KL, -7.877) (tKL, -10.696)
};
\addlegendentry{$\texttt{TC}_{\text{after}}$ (Violates $\mathrm{TC} \ge 0$)}
\end{axis}
\end{tikzpicture}
\caption{Total Correlation values before ($\texttt{TC}_{\text{before}}$, hollow) and after ICA ($\texttt{TC}_{\text{after}}$, filled) for each estimator. Negative post-ICA TC values indicate estimator bias or variance, not true independence. PSS is the only estimator reporting a stable, low, and nonnegative $\texttt{TC}_{\text{after}}$, aligning with theoretical expectations.}
\label{fig:tc_change_plot}
\end{figure}
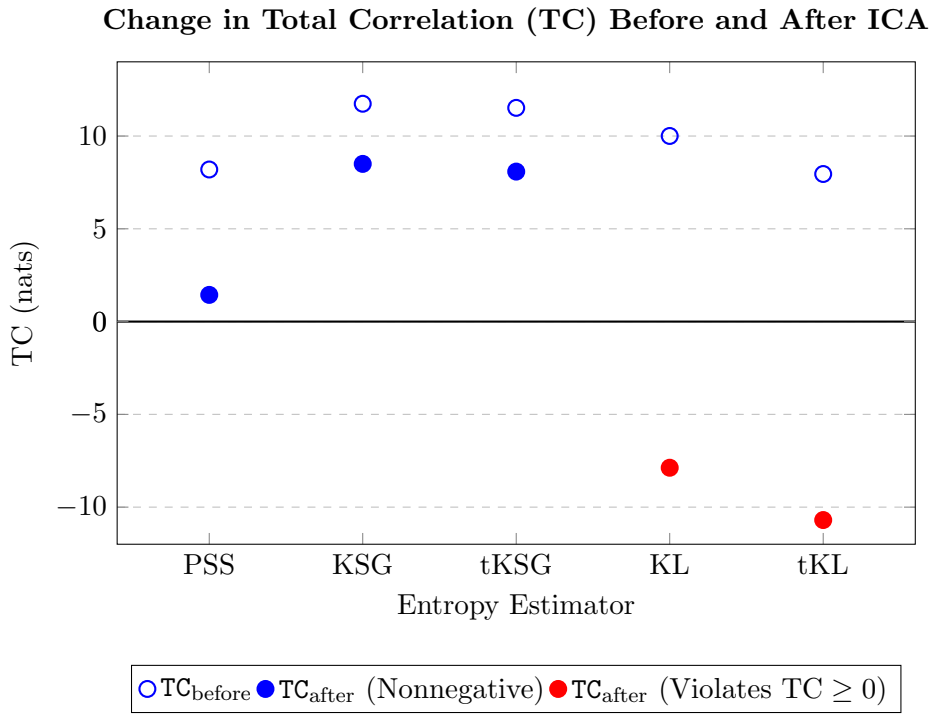

The results are visualized in Figure~\ref{fig:tc_change_plot} and summarized in Table~\ref{tab:results_summary}. These highlight the performance of PSS relative to several widely used $k$-NN based estimators, including KL, KSG, and their truncated variants using normalizing flow (tKL, tKSG).

\begin{table}[h]
\centering
\caption{Summary of ICA performance evaluation by entropy estimator. All TC values are in nats. $\Delta \mathrm{TC} = \texttt{TC}_{\text{before}} - \texttt{TC}_{\text{after}}$.}
\label{tab:results_summary}
\resizebox{\textwidth}{!}{%
\begin{tabular}{l c c c c l}
\toprule
\textbf{Estimator} & \textbf{Parameter} & \textbf{$\texttt{TC}_{\text{after}}$} & \textbf{$\Delta$TC} & \textbf{Time (s)} & \textbf{Verdict} \\
\midrule
\textbf{PSS}  & $\ell^* = 12$ & \textbf{1.439}  & \textbf{6.758}  & \textbf{1.11} & \textbf{Reliable and Sensitive} \\
KSG   & $k^* = 1$ & 8.499  & 3.241  & 1.91 & Moderate Sensitivity \\
tKSG  & $k^* = 1$ & 8.084  & 3.434  & 2.99 & Moderate Sensitivity \\
\midrule
KL    & $k^* = 1$ & $-7.877$  & 17.881  & 1.14 & Violates $\mathrm{TC} \ge 0$ (Estimation Error) \\
tKL   & $k^* = 1$ & $-10.697$ & 18.651  & 0.34 & Violates $\mathrm{TC} \ge 0$ (Estimation Error) \\
\bottomrule
\end{tabular}}
\end{table}

\newpage
\subsubsection{Application to Feature Selection}

In this experiment, we evaluate the practical utility of PSS for feature selection in high-dimensional data. Unlike static entropy estimation, feature selection requires iteratively evaluating mutual information for hundreds of candidate subsets, demanding both computational efficiency and ranking stability. Neural MI estimators, while powerful, require repeated neural-network
optimization and are therefore computationally expensive in feature-selection
settings \citep{Sreejith}, so we restrict our comparison to
nonparametric estimators.

\paragraph{1. Dataset Description.}
We utilized the UCI \texttt{Appliances Energy Prediction} dataset  \\ \citep{appliances_energy_prediction_374}, a widely used benchmark for high-dimensional regression and classification. The dataset consists of $N=19{,}735$ samples and $d=26$ continuous features, including sensors for temperature (T1--T9), humidity (RH\_1--RH\_9), and weather conditions. The features exhibit strong multicollinearity (e.g., correlations between indoor sensors), posing a challenge for distance-based estimators. To frame this as a classification task for mutual information estimation $I(\mathbf{S}; Y)$, we discretized the target variable (energy consumption) into binary classes (High/Low) based on its median.

\paragraph{2. Experimental Setup.}
We performed \emph{Greedy Forward Selection} to rank the features. Starting with an empty set, at each step $t$, we selected the feature $X_j$ that maximized the estimated mutual information with the target $Y$, conditioned on the set $\mathbf{S}_{t-1}$:
\begin{equation*}
    X^* = \operatorname*{argmax}_{X_j \notin \mathbf{S}_{t-1}} \hat{I}(\mathbf{S}_{t-1} \cup \{X_j\}; Y) \, .
\end{equation*}
We compared PSS against the $k$-NN based KL estimator. For PSS, the optimal partition parameter selected via cross-validation was $\ell^*=2$, which effectively adapts to the sparsity of the high-dimensional joint space. Note that neural estimation methods (e.g., MINE, Normalizing Flows) were excluded from this experiment due to their prohibitive computational cost and stochastic instability in iterative search processes.
For evaluation, the selected feature subsets were fed into a Support Vector Machine (SVM) classifier with a Radial Basis Function (RBF) kernel. The classification accuracy was measured on a held-out test set (30\%).

\begin{figure}[h]
    \centering
    % (a) Accuracy Plot
    \begin{minipage}{0.43\textwidth}
        \centering
        \includegraphics[width=\linewidth]{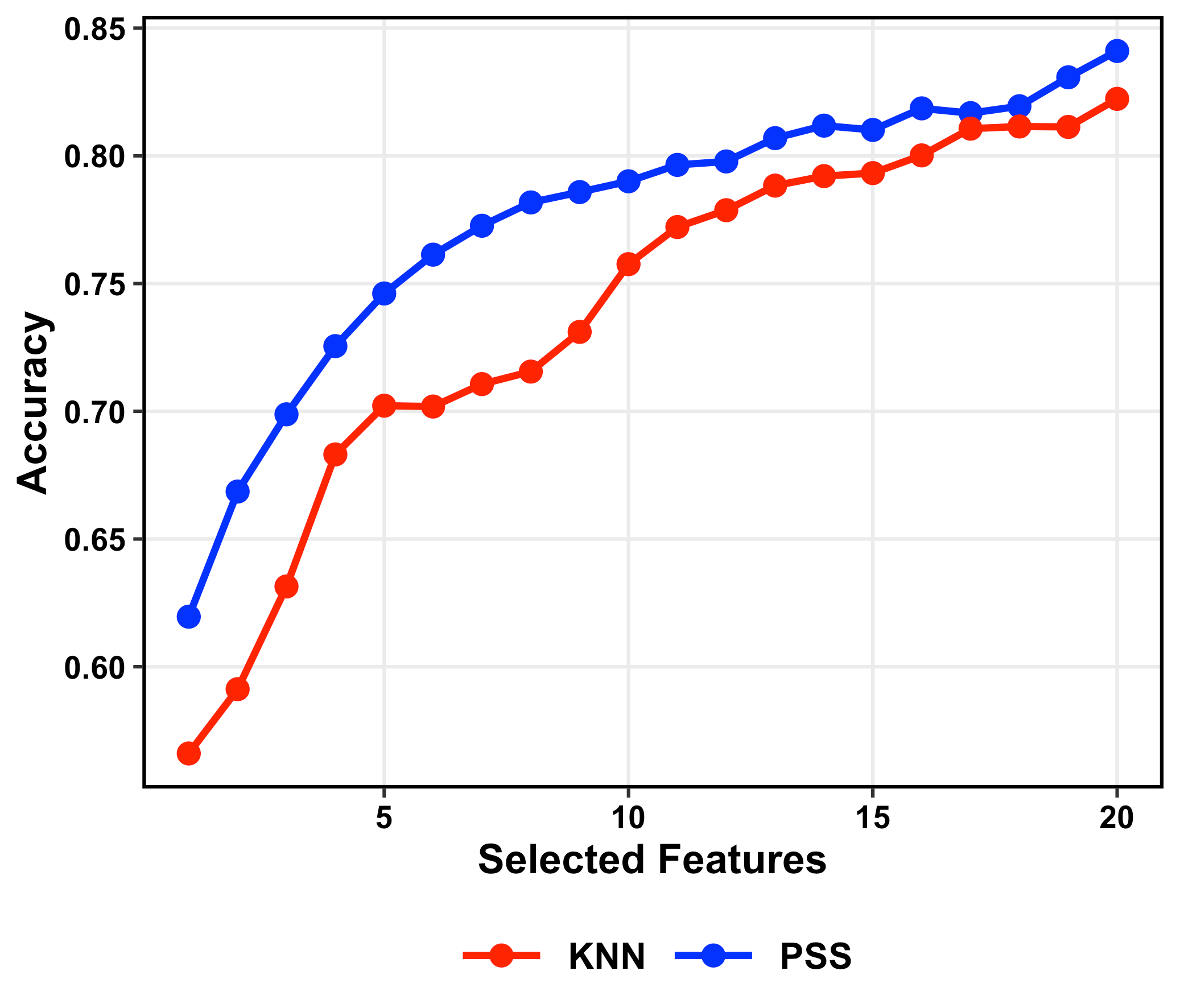} % R에서 저장한 파일명
        \centerline{(a) Classification Accuracy (SVM)}
    \end{minipage}
    \hfill
    % (b) MI Stability Plot
    \begin{minipage}{0.43\textwidth}
        \centering
        \includegraphics[width=\linewidth]{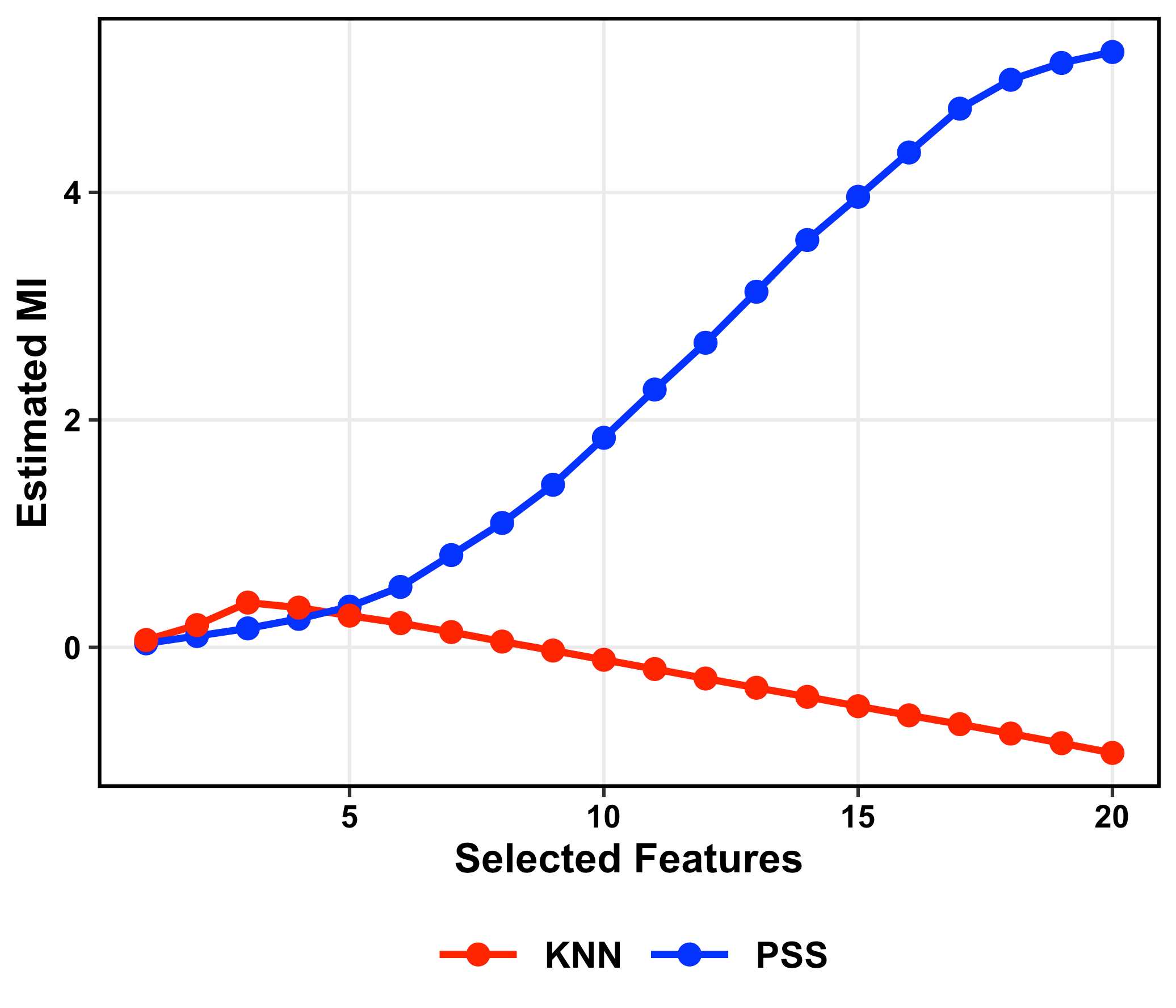} % R에서 저장한 파일명
        \centerline{(b) Estimated MI Stability}
    \end{minipage}
    \caption{Feature selection results on the UCI Energy dataset ($d=26$). \textbf{(a)} SVM classification accuracy. PSS (blue) yields higher accuracy than kNN (red). \textbf{(b)} Estimated Mutual Information. Unlike kNN, which degrades into negative values, PSS maintains theoretical consistency (monotonic increase).}
    \label{fig:feature_selection}
\end{figure}

\paragraph{3. Results and Analysis.}
Figure~\ref{fig:feature_selection} presents the results. 
Figure~\ref{fig:feature_selection}(a) shows the classification accuracy as features are added. PSS (blue) consistently identifies more informative features in the early to middle stages compared to kNN (red), achieving higher accuracy with fewer variables. 

More importantly, Figure~\ref{fig:feature_selection}(b) reveals the theoretical behavior of the estimators. As the dimensionality of the selected subset increases, the estimated MI should fundamentally be non-decreasing (Data Processing Inequality). PSS adheres to this property, showing a stable monotonic increase. In contrast, the kNN estimator collapses after approximately 5 features, producing decreasing and even negative MI values. This confirms that distance-based estimators suffer from severe bias in high-dimensional subspaces ($d > 5$) due to the scarcity of samples, whereas the partitioning strategy of PSS maintains robustness and theoretical consistency.

\newpage
\section{Discussion}
We introduced the Partitioned Sample-Spacing (PSS) estimator, a nonparametric method for joint entropy estimation. By applying univariate spacing methods locally within partitions, PSS achieves almost sure and $L^1$ consistency while avoiding the $O(N^2)$ distance computations of kNN methods. The resulting $O(Nd)$ complexity makes PSS competitive with neural estimators in runtime, without requiring training.

Our experiments show PSS performs well across diverse settings. In synthetic benchmarks on both Gaussian and skewed distributions, PSS consistently achieves lower MSE than kNN methods while matching the accuracy of state-of-the-art normalizing flow variants (UM-tKL, UM-tKSG). Notably, PSS demonstrates remarkable robustness under strong correlation, where its adaptive partitioning strategy effectively captures the concentrated probability mass that typically degrades the performance of global distance-based estimators. We used oracle tuning in these controlled experiments to establish PSS's theoretical potential and enable fair comparison with optimally-tuned baselines.

Real-world applications reveal a critical advantage: PSS maintains theoretical validity where kNN estimators systematically fail. This aligns with the observation of \citep{brown12a} that accurate MI estimation becomes unreliable beyond moderate dimensions, whereas PSS remains stable in exactly these regimes. Using cross-validation for hyperparameter selection, PSS produced non-negative total correlation values in ICA on EEG data ($d=14$), while kNN yielded physically impossible negative values—a fundamental violation of information theory. In feature selection on energy prediction data ($d=26$), PSS preserved the monotonic increase of mutual information as features were added, consistent with the data processing inequality, while kNN exhibited decreasing and negative MI beyond 5 dimensions. These failures of kNN are not mere performance degradation but theoretical collapse in high-dimensional regimes. The features selected by PSS also led to higher downstream classification accuracy, confirming that theoretical soundness translates to practical utility.

PSS has limitations common to grid-based methods. The number of partitions grows as $\ell^d$, making the method impractical when $d \gg 50$ and samples become too sparse. Our experiments show reliable performance up to $d=40$ with moderate sample sizes, which covers many practical applications. Additionally, as noted in the correlation experiments, extreme dependence (e.g., $\rho \to 1$) leading to manifold collapse remains a challenge for full-dimensional density estimators. From a theoretical perspective, we proved consistency but did not derive convergence rates—characterizing the optimal $\ell(n)$ and bias-variance tradeoff remains open. Other directions include adaptive partition size strategies, extension to mixed discrete-continuous variables, and improved scaling to higher dimensions.

Overall, PSS provides a simple, consistent, and training-free approach to multivariate entropy estimation. Its combination of theoretical guarantees, computational efficiency, and robustness to non-Gaussian distributions makes it useful for applications in feature selection, causal discovery, and representation learning where accurate entropy estimates are needed.

\appendix
\section{Proof of Lemma 1}\label{proof:1}
We first claim $X_{(i+m)}-X_{(i)}\xrightarrow{a.s.}0$. Let $\hat{F}_n(x)$ be an empirical cdf of $F(X)$. Then, $\hat{F}_n(X_{(i+m)})=\frac{1}{n}\sum_{j=1}^n I(X_{j}\le X_{(i+m)})=\frac{i+m}{n}$. Note that $\sup\limits_{x}|F(x)-\hat{F}_n(x)|\xrightarrow{a.s.}0$ by Glivenko-Cantelli theorem. Then, $ |F(X_{(i)})-F(X_{(i+m)})|$ is bounded by:
\begin{align*}
    |F(X_{(i)})-\hat{F}_n(X_{(i)})|+|\hat{F}_n(X_{(i)})-\hat{F}_n(X_{(i+m)})|+|\hat{F}_n(X_{(i+m)})-F(X_{(i+m)})|
\end{align*}
The first and third term goes to 0 almost surely by Glivenko-Cantelli theorem. The second term goes to 0 deterministically as $|\hat{F}_n(X_{(i)})-\hat{F}_n(X_{(i+m)})|=\frac{m}{n}$ which goes to 0 by the assumption. Hence,
\begin{equation*}
    F(X_{(i+m)})-F(X_{(i)})\xrightarrow{a.s.}0
\end{equation*}
Since $F$ is continuous, $X_{(i+m)}-X_{(i)}\xrightarrow{a.s.}0$. Similarly, $X_{(i-m)}-X_{(i)}\xrightarrow{a.s.}0$. \\Since $\xi_i=\frac{1}{2m}\sum\limits_{m'=-m}^{m-1}X_{(i+m')}$, 
\begin{align*}
    |\xi_{i}-X_{(i)}|\ = |\frac{1}{2m}\sum\limits_{m'=-m}^{m-1}(X_{(i+m')}-X_{(i)})|&\le \frac{1}{2m}\sum\limits_{m'=-m}^{m-1}|X_{(i+m')}-X_{(i)}|\\
    &\le \frac{1}{2m}\sum\limits_{m'=-m}^{m-1}|X_{(i+m)}-X_{(i-m)}|\\
    &=|X_{(i+m)}-X_{(i-m)}|
\end{align*}
As $X_{(i+m)}-X_{(i-m)}\xrightarrow{a.s.}0$, $\xi_i-X_{(i)}\xrightarrow{a.s.}0$. Similarly, as $\xi_{i+1}=\frac{1}{2m}\sum\limits_{m'=-m+1}^{m}X_{(i+m')}$, $\xi_{i+1}-X_{(i)}\xrightarrow{a.s.}0$.

\section{Proof of Theorem 2}\label{proof:2}

The estimator for $x \in (\xi_i, \xi_{i+1}]$ can be written as the product of two terms:
\[
    \hat{f}_n(x) = \frac{2m/n}{L_n} = \left( \frac{2m/n}{P_n} \right) \cdot \left( \frac{P_n}{L_n} \right),
\]
where $L_n = X_{(i+m)} - X_{(i-m)}$ and $P_n = F(X_{(i+m)}) - F(X_{(i-m)})$. The proof proceeds by showing that the first term in the product converges to 1 almost surely, while the second term converges to $f(x)$ almost surely.

\textbf{Convergence of $P_n$}

Let $U_k = F(X_{(k)})$. Since $F$ is continuous, $U_k$ are the order statistics from a Uniform(0,1) distribution.
Thus, $P_n = U_{(i+m)} - U_{(i-m)}$.
The random variable $P_n = U_{(i+m)} - U_{(i-m)}$ follows a Beta distribution, specifically $P_n \sim \mathrm{Beta}(2m, n-2m+1)$ \cite{david2004order}.

This is a standard result for differences of uniform order statistics. If $U_{(1)} < \dots < U_{(n)}$ are order statistics from $U(0,1)$, then for $j>k$, $U_{(j)} - U_{(k)} \sim \mathrm{Beta}(j-k, n-(j-k)+1)$.
In our case, we have $U_{(i+m)} - U_{(i-m)}$. Let $j' = i+m$ and $k' = i-m$.
The first parameter is $\alpha = j' - k' = (i+m) - (i-m) = 2m$.
The second parameter is $\beta = n - (j'-k') + 1 = n - 2m + 1$.
Thus, $P_n \sim \mathrm{Beta}(2m, n-2m+1)$. Let $\mu_n$ be the mean of $P_n$, $\mathbb{E}[P_n] = \frac{2m}{2m + (n-2m+1)} = \frac{2m}{n+1}$. 

We will now prove that the ratio $\frac{P_n}{\mu_n}$ converges to 1 almost surely. To do so, we will use the Borel-Cantelli Lemma, which requires us to show that the sum of the probabilities of deviation, $\sum_{n=1}^\infty \mathbb{P}\left(\left| \frac{P_n}{\mu_n} - 1 \right| > \epsilon\right)$, is finite for any arbitrary $\epsilon > 0$.
We employ a standard concentration inequality for the Beta distribution, which is sharper than the general Hoeffding bound. For a random variable $X \sim \text{Beta}(\alpha, \beta)$ with mean $\mu$, a bound can be derived from the Chernoff method \cite{Massart}. One such bound is given by:
\begin{equation*}
\mathbb{P}(|X - \mu| \ge t) \le 2 \exp\left(-\frac{(\alpha+\beta)t^2}{2\mu(1-\mu) + \frac{2}{3}t}\right)
\end{equation*}

For our problem, we apply this with $X=P_n$, $\alpha+\beta = n+1$, $\mu=\mu_n$, and $t=t_n=\mu_n\epsilon$. The probability of the event $A_n(\epsilon) = \{|P_n - \mu_n| > \mu_n\epsilon\}$ is thus bounded by:
\begin{align*}
    \mathbb{P}(A_n(\epsilon)) \le 2 \exp\left(-\frac{(n+1)(\mu_n\epsilon)^2}{2\mu_n(1-\mu_n) + \frac{2}{3}\mu_n\epsilon}\right) &= 2 \exp\left(-\frac{2m\epsilon^2}{2(1-\mu_n) + \frac{2}{3}\epsilon}\right)\\
    &\le 2\exp\left( -\frac{\epsilon^2\sqrt{n}}{1 + \frac{1}{3}\epsilon} \right)
\end{align*}
 
Let $K = \frac{\epsilon^2}{1 + 1/3\epsilon}$. We use the Limit Comparison Test with the convergent p-series $\sum_{n=1}^\infty \frac{1}{n^2}$. We evaluate the limit of the ratio:
$$ \lim_{n\to\infty} \frac{2\exp(-K\sqrt{n})}{1/n^2} = 2 \lim_{n\to\infty} \frac{n^2}{e^{K\sqrt{n}}} = 0 $$
Since the limit is 0 and $\sum \frac{1}{n^2}$ converges, the series $\sum_{n=1}^\infty 2\exp(-K\sqrt{n})$ also converges.

Hence,
\begin{align*}
    \sum_{n=1}^\infty \mathbb{P}\left(\left| \frac{P_n}{\mu_n} - 1 \right| > \epsilon\right)
    =&\sum_{n=1}^\infty\mathbb{P}\left(\left| P_n - \mu_n \right| > \mu_n\epsilon\right)\\=&\sum_{n=1}^\infty\mathbb{P}(A_n(\epsilon))
    \le\sum_{n=1}^\infty 2\exp(-K\sqrt{n})<\infty
\end{align*}

By Borel-Cantelli Lemma, $ \frac{P_n}{\mu_n} \xrightarrow{a.s.} 1 $. Its reciprocal $\frac{\mu_n}{P_n}$ also converges to 1 almost surely by the continuous mapping theorem. Finally, since the deterministic sequence $\frac{n+1}{n}$ converges to 1, the almost sure convergence of $\frac{\mu_n}{P_n}$ to 1 directly implies the convergence of their product:
\[ \frac{2m/n}{P_n} = \left(\frac{\mu_n}{P_n}\right) \cdot \left(\frac{n+1}{n}\right) \xrightarrow{a.s.} 1 \]
\\
\textbf{Convergence of $\frac{P_n}{L_n}$ to $f(x)$}

We have $\frac{P_n}{L_n} = \frac{F(X_{(i+m)}) - F(X_{(i-m)})}{X_{(i+m)} - X_{(i-m)}}$.
Then by the Mean Value Theorem, there exists $X'_i \in (X_{(i-m)}, X_{(i+m)})$ such that
$$ \frac{F(X_{(i+m)}) - F(X_{(i-m)})}{X_{(i+m)} - X_{(i-m)}} = f(X'_i) $$

Since $x \in (\xi_i, \xi_{i+1}]$, and $\xi_{i+1} - \xi_i \xrightarrow{a.s.} 0$ (by Lemma 1), for $x$ to remain in this asymptotically vanishing interval, it must be that $\xi_i \xrightarrow{a.s.} x$.
Given $\xi_i \xrightarrow{a.s.} X_{(i)}$ and $\xi_i \xrightarrow{a.s.} x$, it follows that $X_{(i)} \xrightarrow{a.s.} x$.

Now, from Lemma 1, $X_{(i-m)} \xrightarrow{a.s.} X_{(i)}$ and $X_{(i)} \xrightarrow{a.s.} x$. By transitivity of almost sure convergence, $X_{(i-m)} \xrightarrow{a.s.} x$.
Similarly, since $X_{(i+m)} \xrightarrow{a.s.} X_{(i)}$ and $X_{(i)} \xrightarrow{a.s.} x$, then $X_{(i+m)} \xrightarrow{a.s.} x$.

We have $X_{(i-m)} \le X'_i \le X_{(i+m)}$.
Since $X_{(i-m)} \xrightarrow{a.s.} x$ and $X_{(i+m)} \xrightarrow{a.s.} x$, by the squeeze theorem, $X'_i \xrightarrow{a.s.} x$. Since $f$ is continuous at $x$, by the Continuous Mapping Theorem, $f(X'_i) \xrightarrow{a.s.} f(x)$. Therefore, $\frac{P_n}{L_n} \xrightarrow{a.s.} f(x)$.

We have $\hat{f}_n(x) = \left( \frac{2m/n}{P_n} \right) \cdot \left( \frac{P_n}{L_n} \right)$.
Since $\frac{2m/n}{P_n} \xrightarrow{a.s.} 1$ and $\frac{P_n}{L_n} \xrightarrow{a.s.} f(x)$, their product converges almost surely:
$$ \hat{f}_n(x) \xrightarrow{a.s.} 1 \cdot f(x) = f(x) $$

\section{Proof of Proposition 3}\label{proof:3}
Let the sample space be divided into $\ell^2$ partitions. The partitioned density estimator $\hat{f}_{n,\ell}$ in a cell $(\xi_i^k, \xi_{i+1}^k] \times (\eta_j^k, \eta_{j+1}^k]$ within partition $P_k$ is defined as
\[
\hat{f}_{n,\ell}(x,y) = \frac{n_k}{n} \cdot \frac{4m_k^2}{n_k^2 \Delta x_{i}^{k} \Delta y_{j}^{k}}
\]
where $\Delta x_{i}^{k}=x_{(i+m_{k})}^{k}-x_{(i-m_{k})}^{k}$ and  $\Delta y_{j}^{k}=y_{(j+m_{k})}^{k}-y_{(j-m_{k})}^{k}$. The total integral is the sum of integrals over each partition. The integral over a single partition $P_k$ can be split into contributions from its $(n_k-1)^2$ interior cells and its boundary cells:\\ $\int_{C_k} \hat{f}_{n,\ell} = InteriorMass_k + EdgeMass_k$. For any interior cell, the ratio of the grid spacing to the m-spacing simplifies. Using the definition $\eta_j^k = \frac{1}{2m_k}\sum_{l=j-m_k}^{j+m_k-1} y_{(l)}^k$, the difference $\eta_{j+1}^k - \eta_j^k$ becomes a telescoping sum:
\begin{align*}
    \eta_{j+1}^k - \eta_j^k &= \frac{1}{2m_k} \left( \sum_{l=j-m_k+1}^{j+m_k} y_{(l)}^k - \sum_{l=j-m_k}^{j+m_k-1} y_{(l)}^k \right) = \frac{1}{2m_k}\left(y_{(j+m_k)}^k - y_{(j-m_k)}^k\right) = \frac{\Delta y_j^k}{2m_k}
\end{align*}
Thus, the ratio $\frac{\eta_{j+1}^{k}-\eta_{j}^{k}}{\Delta y_{j}^{k}}$ is exactly $\frac{1}{2m_k}$. Similarly, $\frac{\xi_{i+1}^{k}-\xi_{i}^{k}}{\Delta x_{i}^{k}}=\frac{1}{2m_k}$.  Then, The interior mass is:

\begin{align*}
InteriorMass_k =&\sum_{i=1}^{n_k-1}\sum_{j=1}^{n_k-1}\left( \frac{n_k}{n}\frac{4m_k^2}{n_k^2\Delta x_{i}^{k}\Delta y_{j}^{k}}\right)\times (\xi_{i+1}^{k}-\xi_{i}^{k})(\eta_{j+1}^{k}-\eta_{j}^{k}) \\
=& \frac{n_k}{n}\frac{4m_k^2}{n_k^2}\left(\sum_{i=1}^{n_k-1}\frac{\xi_{i+1}^{k}-\xi_{i}^{k}}{\Delta x_{i}^{k}}\right)\left(\sum_{j=1}^{n_k-1}\frac{\eta_{j+1}^{k}-\eta_{j}^{k}}{\Delta y_{j}^{k}}\right) = \frac{n_k}{n}\left(1-\frac{1}{n_k}\right)^2
\end{align*}

The edge mass depends on boundary terms $E_{x}^{k} := \frac{\xi_{1}^{k}-\xi_{0}^{k}}{\Delta x_{0}^{k}}+\frac{\xi_{n_{k}+1}^{k}-\xi_{n_{k}}^{k}}{\Delta x_{n_{k}}^{k}}$ and \\$E_{y}^{k} := \frac{\eta_{1}^{k}-\eta_{0}^{k}}{\Delta y_{0}^{k}}+\frac{\eta_{n_{k}+1}^{k}-\eta_{n_{k}}^{k}}{\Delta y_{n_{k}}^{k}}$. Using the definitions of $\xi_i^k$ and the boundary conditions on order statistics, we can show $0 \le E_{x}^{k}, E_{y}^{k} \le 1$. The first term of $E_x^k$ is bounded as:
\[
\frac{\xi_{1}^{k}-x_{(1)}^{k}}{\Delta x_{0}^{k}} = \frac{\sum_{j=1}^{m_k}(x_{(j)}^k-x_{(1)}^{k})}{2m_k\Delta x_{0}^{k}} \le \frac{m_k(x_{(m_k)}^k-x_{(1)}^{k})}{2m_k\Delta x_{0}^{k}} = \frac{1}{2}
\]
A similar argument for the second term shows $E_y^k \le 1/2 + 1/2 = 1$. The edge mass denoted as $ EdgeMass_k$ is the sum of the remaining terms from the factorization:
\begin{align*}
   &\frac{n_k}{n}\frac{4m_k^2}{n_k^2} \left( \sum_{i=1}^{n_{k}-1}\frac{\xi_{i+1}^{k}-\xi_{i}^{k}}{\Delta x_{i}^{k}} \right)E_y^k+\frac{n_k}{n}\frac{4m_k^2}{n_k^2} \left( \sum_{j=1}^{n_{k}-1}\frac{\eta_{j+1}^{k}-\eta_{j}^{k}}{\Delta y_{j}^{k}} \right)E_x^k+\frac{n_k}{n}\frac{4m_k^2}{n_k^2}E_x^kE_y^k\\ 
    =& \frac{n_k}{n}\frac{4m_k^2}{n_k^2}\left[ \left(\frac{n_k-1}{2m_k}\right)E_y^k + \left(\frac{n_k-1}{2m_k}\right)E_x^k + E_x^k E_y^k \right] \\
    \le& \frac{n_k}{n}\left[ \frac{2m_k(n_k-1)}{n_k^2} + \frac{2m_k(n_k-1)}{n_k^2} + \frac{4m_k^2}{n_k^2} \right] = \frac{n_k}{n}O\left(\frac{m_k}{n_k}\right)
\end{align*}
Summing over all partitions:
\begin{align*}
\int_{\mathbb{R}^{2}}\hat{f}_{n,\ell} &= \sum_{k=1}^{\ell^2} InteriorMass_k + \sum_{k=1}^{\ell^2} EdgeMass_k \\
&= \sum_{k=1}^{\ell^2} \frac{n_k}{n}\left(1-\frac{1}{n_k}\right)^2 + \sum_{k=1}^{\ell^2} \frac{n_k}{n}O\left(\frac{m_k}{n_k}\right) \\
&= 1 - O\left(\frac{\ell^2}{n}\right) + O\left(\max_k \frac{m_k}{n_k}\right)
\end{align*}
As $n \to \infty$, since $\ell(n)^2 = o(n)$ and $m_k/n_k \to 0$, the total integral converges to 1.

\section{Proof of Lemma 4}\label{proof:4}
We aim to show that $\hat{f}_{X}^k(x) \xrightarrow{a.s.} f_{X|C_k}(x)$ for $x \in (\xi_i^k, \xi_{i+1}^k]$, where $f_{X|C_k}(x)$ is the true conditional probability density function of $X$ given partition $P_k$. The estimator is
$$ \hat{f}_{X}^k(x) = \left( \frac{2m_k/n_k}{P_n^k} \right) \cdot \left( \frac{P_n^k}{L_n^k} \right) $$
where $P_n^k = F_{X|C_k}(X_{(i+m_k)}^k) - F_{X|C_k}(X_{(i-m_k)}^k)$ and $L_n^k = X_{(i+m_k)}^k - X_{(i-m_k)}^k$. The samples $X_{(j)}^k$ are the $n_k$ order statistics from partition $P_k$. This proof relies on the conditions $n_k \to \infty$, $m_k \to \infty$, $m_k/n_k \to 0$ as $n \to \infty$, and continuity of $f_{X|C_k}(x)$ at $x$.

The argument closely follows that of Theorem 1. 

\noindent
\textbf{Convergence of the first term:}
Let $U_{(j)}^k = F_{X|C_k}(X_{(j)}^k)$. These are the order statistics from a Uniform(0,1) distribution, based on the $n_k$ samples in $C_k$.
Thus, \\$P_n^k = U_{(i+m_k)}^k - U_{(i-m_k)}^k \sim \mathrm{Beta}(2m_k, n_k-2m_k+1)$.
Let $\mu_n^k$ be expectation $\mathbb{E}[P_n^k] = \frac{2m_k}{n_k+1}$. Then,
\begin{align*}
    \sum_{n=1}^\infty \mathbb{P}\left(\left| \frac{P_n^k}{\mu_n^k} - 1 \right| > \epsilon\right)=\sum_{n=1}^\infty\mathbb{P}\left(\left| P_n^k - \mu_n^k \right| > \mu_n^k\epsilon\right)\le\sum_{n=1}^\infty 2\exp\left( -\frac{\epsilon^2\sqrt{n_k}}{1 + \frac{1}{3}\epsilon} \right)<\infty
\end{align*}
By Borel-Cantelli Lemma, $ \frac{P_n^k}{\mu_n^k} \xrightarrow{a.s.} 1 $. Its reciprocal $\frac{\mu_n^k}{P_n^k}$ also converges to 1 almost surely by the continuous mapping theorem. Finally, since the deterministic sequence $\frac{n_k+1}{n_k}$ converges to 1, $\frac{2m_k/n_k}{P_n^k}  \xrightarrow{a.s.} 1$

\noindent
\textbf{Convergence of the second term:}
By the Mean Value Theorem (assuming $F_{X|C_k}$ is differentiable with derivative $f_{X|C_k}$):
$$ \frac{P_n^k}{L_n^k} = \frac{F_{X|C_k}(X_{(i+m_k)}^k) - F_{X|C_k}(X_{(i-m_k)}^k)}{X_{(i+m_k)}^k - X_{(i-m_k)}^k} = f_{X|C_k}(X'_{k}) $$
for some $X'_{k} \in (X_{(i-m_k)}^k, X_{(i+m_k)}^k)$.
We need to show $X'_{k} \xrightarrow{a.s.} x$. The interval $(\xi_i^k, \xi_{i+1}^k]$, in which $x$ lies, is defined using averages of order statistics $X_{(j)}^k$ partition $k$, similar to $\xi_i$ in Lemma 1.
Applying the logic of Lemma 1 to the $n_k$ samples within partition $k$:\\
$X_{(i \pm m_k)}^k \xrightarrow{a.s.} X_{(i)}^k$ (since $m_k/n_k \to 0$). Furthermore, for $x \in (\xi_i^k, \xi_{i+1}^k]$, $X_{(i)}^k$ is the relevant order statistic that converges almost surely to $x$.
Since $X'_{k}$ is between $X_{(i-m_k)}^k$ and $X_{(i+m_k)}^k$, both of which converge almost surely to $x$, it follows by the Sandwich Theorem that\\ $X'_{k} \xrightarrow{a.s.} x$.
Given that $f_{X|C_k}$ is continuous at $x$, $f_{X|C_k}(X'_{k}) \xrightarrow{a.s.} f_{X|C_k}(x)$ by the Continuous Mapping Theorem. Finally, the product of the two terms converges almost surely:
$$ \hat{f}_{X}^k(x) = \left( \frac{2m_k/n_k}{P_n^k} \right) \cdot \left( \frac{P_n^k}{L_n^k} \right) \xrightarrow{a.s.} 1 \cdot f_{X|C_k}(x) = f_{X|C_k}(x) $$

\section{Proof of Theorem 5}\label{proof:5}

Let the partition $k$ containing the point $(x,y)$ be defined by the rectangular region \\$[x_{\text{min},k}, x_{\text{max},k}] \times [y_{\text{min},k}, y_{\text{max},k}]$. Let $\Delta x_k = x_{\text{max},k} - x_{\text{min},k}$ be the width and \\$\Delta y_k = y_{\text{max},k} - y_{\text{min},k}$ be the height of this partition, so its area is $A_k = \Delta x_k \Delta y_k$.
The estimator for $f(x,y)$ at the point $(x,y)$ is: $\hat{f}_{n,l}(x,y) = \frac{n_k}{n} \cdot \hat{f}_{X}^k(x) \cdot \hat{f}_{Y}^k(y)$.
We aim to show that $\hat{f}_{n,l}(x,y) \xrightarrow{a.s.} f(x,y)$ as $n,\ell \to \infty$.\\

\noindent
\textbf{Convergence of the Sample Proportion $\frac{n_k}{n}$}

Let $P(C_k) = \iint_k f(u,s) \,du\,ds = \int_{y_{\text{min},k}}^{y_{\text{max},k}} \int_{x_{\text{min},k}}^{x_{\text{max},k}} f(u,s) \,du\,ds$ be the true probability mass of partition $k$. Because the indicator of falling into the partition $k$ is Bernoulli with mean $P(C_k)$ by Kolmogorov's Strong Law of Large Numbers, the sample proportion $\frac{n_k}{n}$ converges almost surely to $P(C_k)$. This relies on the conditions in Lemma 2 that ensure $n_k \to \infty$ for the shrinking partition $k$ containing $(x,y)$.\\

\noindent
\textbf{Consistency of Univariate Density Estimators $\hat{f}_{X}^k(x)$ and $\hat{f}_{Y}^k(y)$}

The true conditional PDF of $X$ at $x$, given $C_k$, is derived as follows:
The conditional PDF $f_{X|C_k}(x)$ is defined from the conditional CDF $F_{X|C_k}(x) = P(X \le x | C_k)$ as:
$$ f_{X|C_k}(x) = \frac{d}{dx} P(X \le x | C_k) = \frac{1}{P(C_k)} \frac{d}{dx} P(X \le x, (X,Y) \in C_k) $$
The joint probability $P(X \le x, (X,Y) \in C_k)$ is:
$$ P(X \le x, (X,Y) \in C_k) = \int_{y_{\text{min},k}}^{y_{\text{max},k}} \left( \int_{x_{\text{min},k}}^{x} f(u,s) \,du \right) ds $$
for $x_{\text{min},k} \le x \le x_{\text{max},k}$. Differentiating with respect to $x$:
\begin{align*} \frac{d}{dx} P(X \le x, (X,Y) \in C_k) &= \int_{y_{\text{min},k}}^{y_{\text{max},k}} \left( \frac{d}{dx} \int_{x_{\text{min},k}}^{x} f(u,s) \,du \right) ds \\ &= \int_{y_{\text{min},k}}^{y_{\text{max},k}} f(x,s) \,ds \end{align*}
by the Fundamental Theorem of Calculus. Thus,
$$ f_{X|C_k}(x) = \frac{\int_{y_{\text{min},k}}^{y_{\text{max},k}} f(x,s) \,ds}{P(C_k)} $$
Similarly, for $Y$:
$$ f_{Y|C_k}(y) = \frac{\int_{x_{\text{min},k}}^{x_{\text{max},k}} f(u,y) \,du}{P(C_k)} $$
The existence and continuity of these conditional PDFs are ensured by the assumption that the joint PDF $f(x,y)$ is continuous and positive, and $P(C_k) > 0$.

According to Lemma 2, under the conditions on $l(n)$, $n_k$, and $m_k$:
$$ \hat{f}_{X}^k(x) \xrightarrow{a.s.} f_{X|C_k}(x) \qquad \hat{f}_{Y}^k(y) \xrightarrow{a.s.} f_{Y|C_k}(y)$$

Since $\frac{n_k}{n} \xrightarrow{a.s.} P(C_k)$, $\hat{f}_{X;k}(x) \xrightarrow{a.s.} f_{X|C_k}(x)$, and $\hat{f}_{Y;k}(y) \xrightarrow{a.s.} f_{Y|C_k}(y)$,  their product also converges almost surely:
$$ \hat{f}_{n,l}(x,y) \xrightarrow{a.s.} P(C_k) \cdot f_{X|C_k}(x) \cdot f_{Y|C_k}(y) $$

\noindent
\textbf{Limit of the Target Expression $P(C_k) \cdot f_{X|C_k}(x) \cdot f_{Y|C_k}(y)$}

Let $Q_k(x,y) = P(C_k) \cdot f_{X|C_k}(x) \cdot f_{Y|C_k}(y)$. Substituting the integral expressions for the conditional densities:
\begin{align*} Q_k(x,y) &= P(C_k) \cdot \left(\frac{\int_{y_{\text{min},k}}^{y_{\text{max},k}} f(x,s) \,ds}{P(C_k)}\right) \cdot \left(\frac{\int_{x_{\text{min},k}}^{x_{\text{max},k}} f(u,y) \,du}{P(C_k)}\right) \\ &= \frac{1}{P(C_k)} \left(\int_{y_{\text{min},k}}^{y_{\text{max},k}} f(x,s) \,ds\right) \left(\int_{x_{\text{min},k}}^{x_{\text{max},k}} f(u,y) \,du\right) \end{align*}
As $n \to \infty$, the conditions (including $l(n) \to \infty$ from Lemma 4) ensure that the partition $k$ containing $(x,y)$ shrinks, so its lengths $\Delta x_k \to 0$ and $\Delta y_k \to 0$. Since $f(u,s)$ is continuous at $(x,y)$:

By the Mean Value Theorem for Integrals, There exists $\tilde{y}_k \in [y_{\text{min},k}, y_{\text{max},k}]$ such that\\ $\int_{y_{\text{min},k}}^{y_{\text{max},k}} f(x,s) \,ds = f(x, \tilde{y}_k) \Delta y_k$. Also, there exists $\tilde{x}_k \in [x_{\text{min},k}, x_{\text{max},k}]$ such that \\$\int_{x_{\text{min},k}}^{x_{\text{max},k}} f(u,y) \,du = f(\tilde{x}_k, y) \Delta x_k$. Lastly, there exists $(\bar{x}_k, \bar{y}_k) \in [x_{\text{min},k}, x_{\text{max},k}] \times [y_{\text{min},k}, y_{\text{max},k}]$ such that $P(C_k) = \iint_k f(u,s) \,du\,ds = f(\bar{x}_k, \bar{y}_k) A_k = f(\bar{x}_k, \bar{y}_k) \Delta x_k \Delta y_k$. Substituting these into the expression for $Q_k(x,y)$:
$$ Q_k(x,y) = \frac{ (f(x, \tilde{y}_k) \Delta y_k) (f(\tilde{x}_k, y) \Delta x_k) }{ f(\bar{x}_k, \bar{y}_k) \Delta x_k \Delta y_k } = \frac{f(x, \tilde{y}_k) f(\tilde{x}_k, y)}{f(\bar{x}_k, \bar{y}_k)} $$
(assuming $f(\bar{x}_k, \bar{y}_k) \neq 0$, which is consistent with $f(x,y) > 0$ and the shrinking partition).

Now, we take the limit as $n \to \infty$, which implies $\Delta x_k \to 0$ and $\Delta y_k \to 0$. Since partition $k$ shrinks to the point $(x,y)$: $\tilde{y}_k \to y, \tilde{x}_k \to x, \bar{x}_k \to x\; \text{and}\; \bar{y}_k \to y$.

Given that $f(u,s)$ is continuous at $(x,y)$:
$$\lim_{n \to \infty} f(x, \tilde{y}_k) = f(x,y),\ \lim_{n \to \infty} f(\tilde{x}_k, y) = f(x,y)\ \text{and}\ \lim_{n \to \infty} f(\bar{x}_k, \bar{y}_k) = f(x,y)$$
Therefore, since $f(x,y) > 0$: $\lim_{n \to \infty} Q_k(x,y) = \frac{f(x,y) \cdot f(x,y)}{f(x,y)} = f(x,y) $

We have established that $\hat{f}_{n,l}(x,y) \xrightarrow{a.s.} Q_k(x,y)$ and that $Q_k(x,y) \to f(x,y)$ as $n \to \infty$ (since $k$ depends on $n$ via $l(n)$ and shrinks to $(x,y)$). Therefore, it follows that:
$$ \hat{f}_{n,l}(x,y) \xrightarrow{a.s.} f(x,y) $$

\section{Proof of Proposition 7}\label{proof:7}

    Let $x_v' = a\,x_v + b,\
y_v' = c\,y_v + d$ for some constant $a,b>0$ and $c,d\in \mathbb{R}$ and note that the grid lines transform as
\(\eta_i^{k'} = a\,\eta_i^k + b\), 
\(\xi_j^{k'} = c\,\xi_j^k + d\).
Then,
\begin{align*}
&H_v(x',y',n,\ell)\\
=&\frac{1}{n}
\sum_{v=1}^n
\sum_{k=1}^{\ell^2}
\sum_{j=1}^{n_k}
\sum_{i=1}^{n_k}
\log\!\Bigl(\tfrac{n\,n_k\,(x'_{(i+m_k)}-x'_{(i-m_k)})
\,(y'_{(j+m_k)}-y'_{(j-m_k)})}{4m_k^2}\Bigr)
\,I(\eta_i^{k'}<x'_v\le\eta_{i+1}^{k'})\,
I(\xi_j^{k'}<y'_v\le\xi_{j+1}^{k'})\\
%% substitute differences
=&\frac{1}{n}
\sum_{k,v,j,i}
\log\!\Bigl(\tfrac{n\,n_k\,
\bigl(a(x_{(i+m_k)}-x_{(i-m_k)})\bigr)\,
\bigl(c(y_{(j+m_k)}-y_{(j-m_k)})\bigr)
}{4m_k^2}\Bigr)
\,I(a\eta_i^{k}<ax_v\le a\eta_{i+1}^{k})\,
I(c\xi_j^{k}<cy_v\le c\xi_{j+1}^{k})\\
=&\frac{1}{n}
\sum_{k,v,j,i}
\Bigl[\log(a\,c)
+\log\!\bigl(\tfrac{n\,n_k\,(x_{(i+m_k)}-x_{(i-m_k)})
\,(y_{(j+m_k)}-y_{(j-m_k)})}{4m_k^2}\bigr)\Bigr]
\,I(\eta_i^{k}<x_v\le\eta_{i+1}^{k})\,
I(\xi_j^{k}<y_v\le\xi_{j+1}^{k})\\
=&\log(a\,c)
\;+\;
\frac{1}{n}
\sum_{k,v,j,i}
\log\!\bigl(\tfrac{n\,n_k\,(x_{(i+m_k)}-x_{(i-m_k)})
\,(y_{(j+m_k)}-y_{(j-m_k)})}{4m_k^2}\bigr)
\,I(\eta_i^{k}<x_v\le\eta_{i+1}^{k})\,
I(\xi_j^{k}<y_v\le\xi_{j+1}^{k})\\
=&\log(a\,c)\;+\;H_v(x,y,n,\ell).
\end{align*}

Since \(\log(a\,c)\) is independent of \(\ell\), it does not affect the location of the minimum.  Hence
\[
\arg\min_\ell H_v(x',y',n,\ell)
=\arg\min_\ell H_v(x,y,n,\ell),
\]
showing that the optimal \(\ell\) is invariant under arbitrary shifts and positive scalings of the data.
Because any change of mean and variance is just a positive scaling followed by a shift, and under any such affine transform $H$ only gains an additive constant, the $\ell$ that minimizes $H$ must be unchanged. Hence, the optimal partitioning parameter $\ell$ is invariant under mean and variance.

\section{Proof of Theorem 8}\label{proof:8}

The proof relies on Vitali's Convergence Theorem via the de la Vallée–Poussin criterion. We have already established from Theorem 5 that $\log\hat{f}_{n,\ell}(X,Y) \xrightarrow{a.s.} \log f(X,Y)$. The remaining critical step is to prove that the sequence of random variables $\{ \log\hat{f}_{n,\ell} \}$ is uniformly integrable. We achieve this by showing it is bounded in $L^{1+\delta}$ for some $\delta \in (0,1]$.

The log-density estimator can be decomposed as:
\[
\log\hat{f}_{n,\ell}(x,y) = \log\left(\frac{n_k}{n}\right) + \log\hat{f}_{X}^k(x) + \log\hat{f}_{Y}^k(y)
\]
where for a point $(x,y)$ in the sub-grid cell $(\xi_i^k, \xi_{i+1}^k] \times (\eta_j^k, \eta_{j+1}^k]$, we define $\hat{f}_{X}^k(x) = S_x T_x$ and $\hat{f}_{Y}^k(y) = S_y T_y$. The terms are:
\begin{align*}
S_x = \frac{2m_k}{n_k\{F_{X|C_k}(x_{(i+m_k)}^k) - F_{X|C_k}(x_{(i-m_k)}^k)\}}, \quad
T_x = \frac{F_{X|C_k}(x_{(i+m_k)}^k) - F_{X|C_k}(x_{(i-m_k)}^k)}{x_{(i+m_k)}^k - x_{(i-m_k)}^k}
\end{align*}
and analogously for $S_y$ and $T_y$. Then, $\log\hat{f}_{n,\ell} = \log(\frac{n_k}{n}) + \log S_x + \log T_x + \log S_y + \log T_y$. By the mean value theorem inside marginal partition interval \(I_y^{(k)}\) and \(I_x^{(k)}\) we have, \\$T_x=\frac{|I_y^{(k)}|}{p_k}\, f(X',y_k^*),  T_y=\frac{|I_x^{(k)}|}{p_k}\,f(x_k^*,Y')$
for some \(X',Y',y_k^*,x_k^*\in C_k\) where \\$p_k=\Prob\{X\in C_k\}$. Using $\log(n_k/(n p_k))+\log p_k=\log(n_k/n)$ and $|C_k|=|I_x^{(k)}||I_y^{(k)}|$ we obtain,
\begin{equation*}
\log\hat f_{n,\ell}(X,Y)
= \log\frac{n_k}{n p_k}
+ \log S_x+\log S_y
+ \log f(X',y_k^*)+\log f(x_k^*,Y')
+ \log\frac{|C_k|}{p_k}
\end{equation*}

By Hölder's inequality, for any real numbers \( x_1, \dots, x_m \) and any \( \varepsilon > 0 \) we have\\ $\left| \sum_{j=1}^m x_j \right|^{1+\varepsilon}
\le m^{\varepsilon} \sum_{j=1}^m |x_j|^{1+\varepsilon}$. Thus,
\begin{align*}
\E[|\log\hat{f}_{n,\ell}(x,y)|^{1+\delta}] 
&\le 6^{\delta} \big( 
    \E[|\log\tfrac{n_k}{np_k}|^{1+\delta}] 
    + \E[|\log S_x|^{1+\delta}] 
    + \E[|\log f(X',y_k^*)|^{1+\delta}] \\
&\quad
    + \E[|\log S_y|^{1+\delta}] 
    + \E[|\log f(x_k^*,Y')|^{1+\delta}]
    + \E[|\log \frac{|C_k|}{p_k}|^{1+\delta}]
\big)
\end{align*}

The proof reduces to showing that each expectation term is uniformly bounded in $n$.\\

\noindent
\textbf{Finiteness of $\E[|\log(n_k/np_k)|^{1+\delta}]$}\\
Let \(U_k:=\frac{n_k-np_k}{np_k}\) so that \(\frac{n_k}{n p_k}=1+U_k\). Also, set an event $\mathcal{A}:=\{|U_k|\le\frac{1}{2}\}$. Then, $\Prob(\mathcal{A}^c)\le 2\exp\{-cnp_k\}$ by Chernoff bound. Now, by Lyapunov Inequality,
\[
\mathbb{E}\big[|U_k|^{1+\delta}\mid C_k\big]\ \le\ \big(\mathbb{E}[U_k^2\mid C_k]\big)^{(1+\delta)/2}
=\Big(\frac{p_k(1-p_k)}{n p_k^2}\Big)^{\!(1+\delta)/2}
\ \le\ C (n p_k)^{-(1+\delta)/2}.
\]
Similarly \(\mathbb{E}[|U_k|^{2(1+\delta)}\mid C_k]\le C (n p_k)^{-(1+\delta)}\). Using the inequality for $|u|\le1/2$, \\ \(|\log(1+u)|\le C(|u|+u^2)\), we can bound $\E[|\log{\frac{n_k}{n}}|^{1+\delta}]$ as below:
\begin{align*}
    \E[|\log{\frac{n_k}{np_k}}|^{1+\delta}]&\le\E[|\log{(1+U_k)}|^{1+\delta}\mid\mathcal{A}]+\E[|\log{(1+U_k)}|^{1+\delta}\mid\mathcal{A}^c]]\\
    &\le2^\delta\E[(|U_k|+U_k^2))^{1+\delta}\mid\mathcal{A}]+2^\delta|\log{np_k}|^{1+\delta}\cdot\Prob(\mathcal{A}^c)\\
    &\le C\sum_k p_k(np_k)^{-(1+\delta)/2}+2^{1+\delta}|\log{np_k}|^{(1+\delta)}e^{-cnp_k}\\
    &\le Cn^{-(1+\delta)/2}\ell^{(1+\delta)}\Big(\sum_k p_k\Big)^{1-(1+\delta)/2}+ \cdots\\
    &\le C(\ell^2/n)^{(1+\delta)/2}+2^{1+\delta}|\log{np_k}|^{(1+\delta)}e^{-cnp_k}
\end{align*}
As $\ell^2=o(n)$, the first term goes to 0. Also the second term goes to 0. Thus, \\$\sup_n \E[|\log(n_k/np_k)|^{1+\delta}]<\infty$.\\

\noindent
\textbf{Finiteness of $\E[|\log S_x|^{1+\delta}]$ and $\E[|\log S_y|^{1+\delta}]$.}\\
We prove this for the $x$ dimension; the proof for the $y$ dimension is identical.\\ Let $P_x^k = F_{X|C_k}(x_{(i+m_k)}^k) - F_{X|C_k}(x_{(i-m_k)}^k)$. This random variable follows a Beta distribution, $P_x^k \sim \mathrm{Beta}(2m_k, n_k-2m_k+1)$. The conditional mean and variance of its logarithm are given in \cite{Johnson}:
\begin{align*}
    \E[\log P_x^k \mid n_k] = \psi(2m_k) - \psi(n_k + 1), \quad
    \Var[\log P_x^k \mid n_k] = \psi_1(2m_k) - \psi_1(n_k + 1)
\end{align*}
where $\psi(\cdot)$ and $\psi_1(\cdot)$ are the digamma and trigamma functions.

Using the standard asymptotic expansions from \cite{NIST}, $\psi(z) = \log{z} - \frac{1}{2z} + O(z^{-2})$ and $\log(1 + \frac{1}{z}) = \frac{1}{z} - \frac{1}{2z^2} + O(z^{-3})$, we obtain
\begin{align*}
\E[\log P_x^k \mid n_k] = \log\left(\frac{2m_k}{n_k}\right) + O\left(\frac{1}{m_k}\right),\quad
\mathrm{Var}[\log{P_x^k}\mid n_k]
= \tfrac{1}{2m_k}+O\!\Big(\tfrac{1}{m_k^2}\Big)
\end{align*}
Now we bound the conditional expectation of $|\log S_x|^{1+\delta}$:
\begin{align*}
    &\E[|\log S_x|^{1+\delta} \mid n_k]\\
    &\le 2^\delta \left| \log\left(\frac{2m_k}{n_k}\right) - \E[\log P_x^k \mid n_k] \right|^{1+\delta} + 2^\delta \E\left[ \left| \log P_x^k - \E[\log P_x^k \mid n_k] \right|^{1+\delta} \mid n_k \right] \\
    &\le 2^\delta \left| O\left( \frac{1}{m_k} \right) \right|^{1+\delta} + 2^\delta \left(\Var\left( \log P_x^k \mid n_k \right)\right)^{\frac{1+\delta}{2}} \\
    &\le O\left( \frac{1}{m_k^{1+\delta}} \right) + 2^\delta \left(\tfrac{1}{2m_k}+O\!\Big(\tfrac{1}{m_k^2}\Big)\right)^{\frac{1+\delta}{2}} \le C n_k^{-(1+\delta)/4},
\end{align*}
using $m_k \asymp \sqrt{n_k}$ and Lyapunov's Inequality. Now set an event $A_k:=\{n_k\ \ge\ \tfrac12\,n p_k\}$. By a Chernoff bound, \(\mathbb{P}(A_k^c\mid C_k)\le \exp(-c\,n p_k)\). On \(A_k\), \(n_k^{-(1+\delta)/4}\le (2/(n p_k))^{r/4}\). Therefore,
\[
\mathbb{E}\big[|\log S_x|^{(1+\delta)}\big]
=\sum_{k} p_k\,\mathbb{E}\big[|\log S_x|^{(1+\delta)}\mid n_k\big]
\le C\sum_k p_k(n p_k)^{-(1+\delta)/4} + C\sum_k p_k e^{-c\,n p_k}.
\]
Using Jensen inequality and $\ell^2=o(n)$, we show the first term goes to 0.
\[ n^{-(1+\delta)/4}\sum_k p_k^{1-(1+\delta)/4}
 \le n^{-(1+\delta)/4}(\ell^2)^{(1+\delta)/4}\Big(\sum_{k=1}^{\ell^2}p_k\Big)^{1-(1+\delta)/4}
 \le (\ell^2/n)^{(1+\delta)/4}
\]
Also $\sum_k p_k e^{-c\,n p_k}\le 1$ for all $n$ (and $\to0$ if $n\min_k p_k\to\infty$). Hence, \(\sup_n \mathbb{E}|\log S_x|^{1+\delta}<\infty\).
The same argument applies to \(\log S_y\).\\

\noindent
\textbf{Finiteness of $\E[|\log f(X',y_k^*)|^{1+\delta}]$, $\E[|\log f(x_k^*,Y')|^{1+\delta}]$, and $\E\big[\big|\log(|C_k|/p_k)\big|^{1+\delta}\big]$}\\
Let $g(u):=|\log f(u)|^{1+\delta}$ with $\int g f<\infty$.
Define the upper step function\\ $s_n(u):=\sum_k\big(\sup_{v\in C_k} g(v)\big)\mathbf 1_{C_k}(u)$.
Since the mesh $\to0$ and $f$ is continuous with $f>0$ a.e., $s_n\downarrow g$ pointwise and
$\int s_n f\downarrow\int g f<\infty$.
Because $(X',y_k^*)\in C_k$ and \\$\bar f_k:=|C_k|^{-1}\!\int_{C_k} f\in[\inf_{C_k} f ,\sup_{C_k} f]$ we have, 
$$|\log f(X',y_k^*)|^{r}\le \sup_{u\in C_k} g(u),\qquad
\Big|\log\frac{|C_k|}{p_k}\Big|^{r}=\big|\log\bar f_k\big|^{r}\le \sup_{u\in C_k} g(u).$$
Hence, $\sup_n$ of each is finite; the same for $\E|\log f(x_k^*,Y')|^{1+\delta}$.

Since each  term in the inequality is uniformly bounded in $L^{1+\delta}$, their sum is also uniformly bounded. Therefore, we have the desired uniform moment control:\\$\sup_{n} \E\left[\,|\log \hat{f}_{n,\ell}(X,Y)|^{1+\delta} \,\right] < \infty$.
Since $1+\delta > 1$, de la Vallée-Poussin's criterion implies that the sequence $\{\log\hat{f}_{n,\ell}\}$ is uniformly integrable. This, combined with the almost sure convergence $\log\hat{f}_{n,\ell} \xrightarrow{a.s.} \log f$, allows us to apply Vitali's Convergence Theorem to conclude that $\lim_{n\to\infty} \E[\log\hat{f}_{n,\ell}] = \E[\log f]$.

This proves that the bias term $B_n := \E[\log\hat{f}_{n,\ell}] - \E[\log f]$ converges to 0. The variance term $A_n := -\frac{1}{n}\sum (\log\hat{f}_{n,\ell}(X_i,Y_i) - \E[\log\hat{f}_{n,\ell}])$ converges to 0 almost surely by the Strong Law of Large Numbers. Thus, $\hat{H}_{n,\ell} \xrightarrow{a.s.} H(f)$.

To prove $L^1$ convergence, we must show the sequence $\{\hat{H}_{n,\ell}\}$ is uniformly integrable. By Jensen's inequality for the convex function $\phi(z)=|z|^{1+\delta}$:
\[
\E\bigl[|\hat{H}_{n,\ell}|^{1+\delta}\bigr]\le \frac{1}{n}\sum_{i=1}^n \E\bigl[|\log\hat{f}_{n,\ell}(X_i)|^{1+\delta}\bigr] = \E\bigl[|\log\hat{f}_{n,\ell}(X_1, Y_1)|^{1+\delta}\bigr]
\]
Since we showed the right-hand side is uniformly bounded in $n$, we have\\ $\sup_n \E\bigl[|\hat{H}_{n,\ell}|^{1+\delta}\bigr] < \infty$. This implies $\{\hat{H}_{n,\ell}\}$ is uniformly integrable. As $\hat{H}_{n,\ell} \xrightarrow{a.s.} H(f)$, Vitali's Convergence Theorem again applies, yielding the final result: $\hat{H}_{n,\ell} \xrightarrow{L^1} H(f)$.

\section{Proof Sketch of Proposition 9}\label{proof:9}
The argument is a direct extension of the proof of Proposition 3.
We now partition $\mathbb{R}^d$ into $\ell^d$ hyperrectangles and define
the local estimator on each cell $P_k$ by
\[
  \hat f_{n,\ell}(x_1,\dots,x_d)
  = \frac{n_k}{n}
    \prod_{j=1}^d \frac{2 m_k}{n_k \Delta x_{j,i_j}^k},
\]
with marginal $m_k$-spacings $\Delta x_{j,i_j}^k$ defined as before.
As in Proposition 3, the total mass decomposes into interior
and boundary contributions:
\[
  \int_{\mathbb{R}^d} \hat f_{n,\ell}
  = \sum_{k=1}^{\ell^d} \text{InteriorMass}_k
  + \sum_{k=1}^{\ell^d} \text{EdgeMass}_k.
\]
For interior cells, each ratio of grid width satisfies
\(
  (\xi_{j,i_j+1}^k - \xi_{j,i_j}^k)/\Delta x_{j,i_j}^k = 1/(2m_k)
\)
by the same telescoping-sum argument as in the bivariate case, yielding\\
\(
  \text{InteriorMass}_k
  = \frac{n_k}{n} (1 - 1/n_k)^d.
\)
The boundary terms can again be written in terms of edge factors
$E_{j}^k$ for $j=1,\dots,d$, each bounded between $0$ and $1$.
Expanding the product and using $0 \le E_{j}^k \le 1$ shows that
\[
  \text{EdgeMass}_k
  = \frac{n_k}{n} \, O\!\left(\frac{m_k}{n_k}\right),
\]
with constants depending only on $d$.
Summing over $k$ gives
\[
  \int_{\mathbb{R}^d} \hat f_{n,\ell}
  = 1 - O\!\left(\frac{\ell^d}{n}\right)
    + O\!\left(\max_k \frac{m_k}{n_k}\right),
\]
and the stated convergence follows from $\ell(n)^d = o(n)$ and $m_k/n_k\to0$.

\section{Proof Sketch of Theorem 10}\label{proof:10}
The proof follows the same three-step logic as the bivariate proof detailed in Appendix~E. First, by the Strong Law of Large Numbers, the sample proportion $\frac{n_k}{n}$ converges almost surely to the true partition probability. Second, by Lemma~4, each of the $d$ univariate conditional density estimators, $\hat{f}_{X_j;k}(x_j)$, converges almost surely to the true conditional density, $f_{X_j|C_k}(x_j)$.

The final step is to show that the limit of the full estimator, $ \frac{n_k}{n} \prod_{j=1}^{d} \hat{f}_{X_j;k}(x_j)$, converges to $f(\mathbf{x})$. By extending the Mean Value Theorem argument from the bivariate case, the product of the limits can be shown to converge to the true joint density.

This establishes the almost sure convergence, $\hat{f}_{n,l}(\mathbf{x}) \xrightarrow{a.s.} f(\mathbf{x})$. The $L^1$ convergence then follows directly from this result and Proposition~9 by an application of Scheffé's Theorem.

\section{Proof Sketch of Theorem 11}\label{proof:11}
The proof establishes the uniform integrability of $\{\log\hat{f}_{n,\ell}(\mathbf{X})\}$ by proving it is bounded in $L^{1+\delta}$. The log-density estimator decomposes as $\log\hat{f}_{n,\ell}(\mathbf{x}) = \log(n_k/n) + \sum_{j=1}^{d} ( \log S_j + \log T_j )$, where
\begin{align*}
S_j &= \frac{2m_k}{n_k\{F_{X_j|C_k}(x_{j,(a_j+m_k)}^k) - F_{X_j|C_k}(x_{j,(a_j-m_k)}^k)\}}, \\
T_j &= \frac{F_{X_j|C_k}(x_{j,(a_j+m_k)}^k) - F_{X_j|C_k}(x_{j,(a_j-m_k)}^k)}{x_{j,(a_j+m_k)}^k - x_{j,(a_j-m_k)}^k}.
\end{align*}
Applying the MVT along the other $(d-1)$ coordinates, $T_j=\frac{|I_{-j}^{(k)}|}{p_k}f(\tilde X_j)$ for some $\tilde X_j\in C_k$. Taking logs and using $\sum_{j=1}^d\log|I_{-j}^{(k)}|=(d-1)\log|C_k|$ yields the final decomposition:
\[
\log\hat f_{n,\ell}(X)
= \underbrace{\log\tfrac{n_k}{n p_k}}_{\text{(A)}}
+ \underbrace{\sum_{j=1}^d \log S_j}_{\text{(B)}}
+ \underbrace{\sum_{j=1}^d \log f(\tilde X_j)}_{\text{(C)}}
+ \underbrace{(d-1)\log\tfrac{|C_k|}{p_k}}_{\text{(D)}}
\]
Let $r=1+\delta$ then, the $L^r$ norms of these terms are bounded as follows, where $s_n$ is the upper step function of $|\log f|^r$:
\begin{align*}
\text{(A)}\quad\E|\log\tfrac{n_k}{n p_k}|^{r}
&\le C(\tfrac{\ell^d}{n})^{r/2}+C,
&
\text{(C)}\quad \E[\sum_{j=1}^d |\log f(\tilde X_j)|^{r}]
&\le d\int s_n f ,
\\
\text{(B)} \sum_{j=1}^d \E|\log S_j|^{r}
&\le C(\tfrac{\ell^d}{n})^{r/4}+C,
&
\text{(D)}\quad \E|(d-1)\log\tfrac{|C_k|}{p_k}|^{r}
&\le (d-1)^r\!\int s_n f.
\end{align*}
where $\int s_n f \downarrow  d\int |\log f|^r f<\infty$. Combining the bounds gives\\ $\sup_n \E|\log\hat f_{n,\ell}(X)|^{r}<\infty$, hence uniform integrability. Almost sure convergence and UI imply $\E[\log\hat f_{n,\ell}(X)]\to \E[\log f(X)]$ by Vitali's theorem. For $\hat H_{n,\ell}$,
a  bounded-difference argument controls the empirical fluctuation,
and Jensen’s inequality implies $\sup_n \E|\hat H_{n,\ell}|^{r}\le \sup_n \E|\log\hat f_{n,\ell}(X)|^{r}<\infty$,
so $\{\hat H_{n,\ell}\}$ is uniformly integrable. Therefore $\hat H_{n,\ell}\to H(f)$ a.s.\ and in $L^1$.

\section*{Code Availability}
All code, scripts, and experimental resources necessary to reproduce the results in this paper are publicly available at: \texttt{https://github.com/hoo0321-design/pss-estimator}

\nocite{Massart, Williams}
            
\bibliography{ref}    

\end{document}